\documentclass{amsart}

\usepackage[english]{babel}
\usepackage[all]{xy}
\usepackage{latexsym}
\usepackage{eufrak}
\usepackage{amssymb}
\usepackage{fancyhdr}

\newtheorem{thm}{\sc Theorem}
\newtheorem{lem}[thm]{\sc Lemma}
\newtheorem{pro}[thm]{\sc Proposition}
\newtheorem{cor}[thm]{\sc Corollary}
\theoremstyle{definition}

\newtheorem*{rem}{\sc Remark}
\newtheorem*{dei}{\sc Definition}
\newtheorem*{deo}{\sc Proof}
\newtheorem*{ex}{\sc Examples}
\newtheorem*{ex1}{\sc Example}

\newcommand{\C}{\mathcal{C}om}
\newcommand{\A}{\mathcal{A}s}
\newcommand{\D}{\mathcal{D}ias}
\newcommand{\Pe}{\mathcal{P}erm}
\newcommand{\Po}{\mathcal{P}}
\newcommand{\Qo}{\mathcal{Q}}
\newcommand{\ac}{\textrm{!`}}
\newcommand{\Li}{\mathcal{L}ie}
\newcommand{\Pli}{\mathcal{P}relie}
\newcommand{\Dend}{\mathcal{D}end}
\newcommand{\F}{\mathcal{F}}

\newcommand{\cqfd}{\ \hfill \square}
\newcommand{\Sy}{\mathbb{S}}
\newcommand{\Sn}{\Sy_n}

\newcommand{\N}{\mathcal{N}}
\newcommand{\B}{\mathcal{B}}

\newcommand{\Co}{\mathcal{C}}
\newcommand{\Trias}{\mathcal{T}rias}
\newcommand{\Tridend}{\mathcal{T}ri\mathcal{D}end}
\newcommand{\Comtri}{\mathcal{C}om\mathcal{T}rias}
\newcommand{\Postlie}{\mathcal{P}ost\mathcal{L}ie}
\newcommand{\Mag}{\mathcal{M}ag}
\newcommand{\U}{\mathcal{U}}

\title{\bf Homology of generalized partition posets}
\author{Bruno Vallette}

\begin{document}

\begin{abstract}
We define a family of posets of partitions associated to an
operad. We prove that the operad is Koszul if and only if the
posets are Cohen-Macaulay. On the one hand, this characterization
allows us to compute completely the homology of the posets. The
homology groups are isomorphic to the Koszul dual cooperad. On the
other hand, we get new methods for proving that an operad is
Koszul.
\end{abstract}

\maketitle

\section*{Introduction}

The homology of the lattice of partitions of the sets $\{1,\,
\ldots,\, n\}$ has been studied for more than 20 years by many
authors. First, A. Bj\"orner proved in \cite{B} that the only
non-vanishing homology groups are in top dimension. Since the
symmetric groups $\Sn$ act on these posets, the homology groups of
top dimension are $\Sn$-modules. It took several years to
completely compute these representations of $\Sn$. We refer to the
introduction of \cite{Fresse} for a complete survey on the
subject. Actually, these homology groups are given by the linear
dual of the multi-linear part of the free Lie algebra twisted by
the signature representation. In \cite{Fresse}, B. Fresse
explained why such a result : the partition lattices are build
upon the operad $\C$ of commutative algebras and the homology of
the partition lattices is isomorphic to its Koszul dual cooperad
$\Li$ corresponding to Lie algebras.\\

In this article, we make explicit a general relation between
operads and partition type posets. From any operad in the category
of sets, we associate a family of partition type posets. For a
class of algebraic operads coming from set operads, we prove the
equivalence of two homological notions : Koszul operad and
Cohen-Macaulay posets.\\

An operad is an algebraic object that represents the operations
acting on certain types of algebras. For instance, one has an $\A$
operad, a $\C$ operad and a $\Li$ operad for associative algebras,
commutative algebras and Lie algebras. An important issue in
 the operadic theory is to show that an operad is Koszul. In this case,
 the associated algebras have interesting properties.

The homology of posets was studied by many authors (\emph{cf.} D.
Quillen \cite{Q}, A. Bj\"orner, A.M. Garsia and R.P. Stanley
\cite{BGS} for instance). In this framework, there exist many
methods to show that the homology of a poset is concentrated in
top dimension and likewise for each interval. In
this case, the poset is called Cohen-Macaulay.\\

The main purpose of this paper is to show that an operad is Koszul
if and only if the related posets are Cohen-Macaulay
(Theorem~\ref{KoszulCohen}). In this case, we can compute the
homology of the poset. This homology is equal to the Koszul dual
cooperad. Since there exist an explicit formula for this dual, one
can use it to compute the homology groups of the posets in terms
of $\Sn$-modules. On the other hand, we can use the combinatorial
methods of poset theory to show that the posets associated to an
operad is Cohen-Macaulay.  Hence, it provides new methods for
proving that an operad is Koszul. Notice that these methods prove
in the same time that the posets are Cohen-Macaulay over any field
$k$ and over the ring of integers $\mathbb{Z}$. Therefore these
new methods work over any field $k$ and over $\mathbb{Z}$. The
classical methods found in the literature are based on the
acyclicity of the Koszul complex defined by the Koszul dual
cooperad. One advantage of the method using the partition posets
is that we do not need to
compute the Koszul dual cooperad (and its coproduct) of an operad to prove that it is Koszul.\\

We introduce several partition type posets (pointed, ordered, with
block size restriction) associated to operads appearing in the
literature. Since these operads are Koszul, we can compute the
homology of the posets. It is concentrated in top dimension and
given by the Koszul dual cooperad. To handle the case of
multi-pointed partitions poset, we introduce a new operad, called
$\Comtri$ for commutative trialgebras. We describe its Koszul dual
and we show that it is a Koszul operad in the appendix.\\

It is worth mentioning here that in \cite{CV}, F. Chapoton and the
author studied the proprieties of the pointed and multi-pointed
posets and applied the results of this paper to prove that the
operads $\Pe$ and $\Comtri$ are Koszul over
$\mathbb{Z}$.\\

Sections $1$ and $2$ contain respectively a survey of the homology
of posets and of Koszul duality for operads. In Section $3$, we
describe the construction of the partitions posets associated to a
(set) operad. We also prove in this section the main theorem of
this paper, namely Theorem \ref{KoszulCohen}, which claims that
the operad is Koszul if and only if the related posets are
Cohen-Macaulay. Examples and applications are treated in Section
$4$. Finally, we prove in the Appendix that the operads $\Comtri$
and $\Postlie$ are
Koszul. \\

Let $k$ be the ring $\mathbb{Z}$, the field $\mathbb{Q}$,
$\mathbb{F}_p$ or any field of characteristic $0$. Unless
otherwise stated, all $k$-modules are assumed to be projective.

\tableofcontents

\section{Order complex of a poset}

We recall the basic definitions of a poset and the example of the
partition poset. (For more details, we refer the reader to Chapter
3 of \cite{Stanley}.) We define the order complex of a poset and
the notion of Cohen-Macaulay.

\subsection{Poset}

\begin{dei}[Poset]
A \emph{poset} $(\Pi,\, \leqslant)$ is a set $\Pi$ equipped with a
partial order relation, denoted by $\leqslant$.
\end{dei}

The posets considered in the sequel are finite. We denote by
$Min(\Pi)$ and $Max(\Pi)$ the sets of minimal and maximal elements
of $\Pi$. When each set $Min(\Pi)$ and $Max(\Pi)$ has only one
element, the poset is said to be \emph{bounded}. In this case, one
denotes by $\hat{0}$ the element of $Min(\Pi)$ and by
$\hat{1}$ the element of $Max(\Pi)$\\

For $x\leqslant y$ in $\Pi$, we denoted the closed interval
$\{z\in \Pi \ | \ x\leqslant z \leqslant y \}$ by $[x,\, y]$ and
the open interval $\{z\in \Pi \ | \ x< z< y \}$ by $(x,\, y)$. For
any $\alpha\in Min(\Pi)$ and any $\omega \in Max(\Pi)$, the closed
interval $[\alpha,\, \omega]$ is a bounded poset. If $\Pi$ is a
bounded poset, one defines the \emph{proper part} $\bar{\Pi}$ of
$\Pi$ by the open interval $( \hat{0},\,
\hat{1})$.\\

For elements $x<y$, if there exists no $z$ such that $x<z<y$, then
we say that $y$ \emph{covers} $x$. The covering relation is
denoted by $x\prec y$.

\begin{dei}[Chain, maximal chain]
A \emph{chain} $\lambda_0< \lambda_1< \cdots<\lambda_l$ is a
growing sequence of elements of a poset $(\Pi,\, \leqslant)$. Its
\emph{length} is equal to $l$.

A \emph{maximal chain} between $x$ and $y$, is a chain
$x=\lambda_0 \prec \lambda_1 \prec \cdots \prec \lambda_l=y$ which
can not be lengthened. A \emph{maximal chain of $\Pi$} is a
maximal chain between an element of $Min(\Pi)$ and an element of
$Max(\Pi)$.
\end{dei}

A poset is \emph{pure} if for any $x\leqslant y$, the maximal
chains between $x$ and $y$ have same length. If a poset is both
bounded and pure, it is called a \emph{graded} poset.

\subsection{Partition Poset}

The set $\{ 1,\ldots,\,n\}$ is denoted by $[n]$.

\begin{dei}[Partitions of $\lbrack n\rbrack $]
A \emph{partition} of the set $[n]$ is an unordered collection of
subsets $B_1,\ldots,\, B_k$, called \emph{blocks} or
 \emph{components}, which are nonempty, pairwise disjoint, and whose union
gives $[n]$.
\end{dei}

\begin{dei}[Partition poset, $\Pi$]
For any integer $n$, one defines a partial order $\leqslant$, on
the set of partitions of $[n]$, by the \emph{refinement of
partitions}. This partially ordered set is called the
\emph{partition poset} (or partition lattice) and denoted by
$\Pi(n)$.
\end{dei}

For instance, one has $\left\{ \{1,\,3\},\,\{2,\,4\}\right\}
\leqslant \left\{ \{1\},\,\{3\},\,\{2,\,4\} \right\}$. The single
set $\{ \{1,\ldots,$ $n\}\}$ forms the smallest partition of $[n]$
whereas the collection $\{ \{1\},\ldots,\{n\} \}$ forms the
largest partition. Therefore, this poset is graded. Observe that
this definition is dual to the one found in the literature.  \\

The set of partitions of $[n]$ is equipped with a right action of
the symmetric group $\Sn$. Let  $\sigma\, :\, \{1,\ldots,\,
n\}\,\rightarrow\,\{1,\ldots,\, n\}$ be a permutation, the image
of the partition $\{\{i_1^1,\ldots ,\, i^1_{j_1}\},\ldots ,\,
\{i^k_1,\ldots,\, i^k_{j_k} \} \}$ under $\sigma$ is the partition
$\{ \{\sigma(i_1^1),\ldots ,\, \sigma(i^1_{j_1})\},\ldots,$ $ $
$\,\{\sigma(i^k_1),\ldots,\,\sigma(i^k_{j_k}) \} \}$.

\subsection{Order complex}

We consider the set of chains $\lambda_0< \lambda_1<
\cdots<\lambda_l$ of a poset $(\Pi,\, \leqslant)$ such that
$\lambda_0\in Min(\Pi)$ and $\lambda_l\in Max(\Pi)$. This set is
denoted by $\Delta(\Pi)$. More precisely, a chain $\lambda_0<
\lambda_1< \cdots<\lambda_l$ of length $l$ belongs to
$\Delta_l(\Pi)$

\begin{dei}[Order complex, $\Delta(\Pi)$]
The set $\Delta(\Pi)$ can be equipped with face maps. For $0<i<l$,
the face map $d_i$ is given by the omission of $\lambda_i$ in the
sequence $\lambda_0 <\lambda_1 <\cdots <\lambda_l$. If we take the
convention $d_0=d_l=0$, the module $k[\Delta(\Pi)]$ is a
presimplicial module. The induced chain complex on
$k[\Delta(\Pi)]$ is called the \emph{order complex} of $\Pi$.
\end{dei}

\begin{dei}[Homology of a poset $\Pi$, $H_*(\Pi)$]
The \emph{homology of a poset $(\Pi,\, \leqslant)$} is the
homology of the presimplicial set $\Delta(\Pi)$ with coefficients
in $k$. We denote it by $H_*(\Pi)$.
\end{dei}

\begin{rem}
In the literature (\emph{cf.} for instance \cite{BGS}), the
\emph{reduced} homology of a poset is defined in the following
way. One denotes by $\widetilde{\Delta}_l(\Pi)$ the set of chains
$\lambda_0< \lambda_1< \cdots<\lambda_l$, with no restriction on
$\lambda_0$ and $\lambda_l$. The face maps $d_i$ are defined by
the omission of $\lambda_i$ for $0\leqslant i \leqslant l$. By
convention, this complex is augmented by
$\widetilde{\Delta}_{-1}(\Pi)=\{ \emptyset \}$. We have
$k[\widetilde{\Delta}_{-1}(\Pi)]=k$. The associated homology
groups are denoted by $\widetilde{H}_*(\Pi)$. This definition is
convenient when applied to the proper part of a bounded poset.

The relation between the two definitions is given by the following
formula
$$\Delta_l(\Pi)=\bigsqcup_{(\alpha,\, \omega) \in Min(\Pi)\times Max(\Pi)}
\widetilde{\Delta}_{l-2}\left(\,(\alpha, \, \omega) \,\right) ,$$
which induces an isomorphism of presimplicial complexes.
Therefore, we have
$$H_l(\Pi) = \bigoplus_{(\alpha,\, \omega) \in Min(\Pi)\times Max(\Pi)}
\widetilde{H}_{l-2}(\, (\alpha, \, \omega)\, ).$$
\end{rem}

When a poset $(\Pi,\, \leqslant)$ is equipped with an action of a
group $G$, compatible with the partial order $\leqslant$, the
modules $k[\Delta_l(\Pi)]$ are $G$-modules. Since, the chain map
commutes with the action of $G$, the homology groups $H_*(\Pi)$
are $G$-modules too. In the case of the partition poset, the
module $k[\Delta(\Pi(n))]$ is an $\Sn$-presimplicial module.

\subsection{Cohen-Macaulay poset}

\begin{dei}[Cohen-Macaulay poset]
Let $\Pi$ be a graded poset. It is said to be
\emph{Cohen-Macaulay} (over $k$) if the homology of each interval
is concentrated in top dimension. For any $x\leqslant y$, let $d$
be the length of maximal chains between $x$ and $y$, we have
$$H_l([x,\, y]) = \widetilde{H}_{l-2}(\, (x,\, y) \,)=0, $$
for $l\neq d$.
\end{dei}

We refer the reader to the article of A. Bj\"orner, A.M. Garsia
and R.P. Stanley \cite{BGS} for a short survey on Cohen-Macaulay
posets.

\section{Operads and Koszul duality}

In this section, we define the notions related to operads. We
recall the results of the Koszul duality theory for operads that
will be used later in the text.

\subsection{Definition of an operad and examples}
We recall the definition of an operad. We give the examples of
free and quadratic operads.

\subsubsection{\bf Definition} The survey of J.-L. Loday \cite{LodayBourbaki} provides
a good introduction to operads and the book of M. Markl, S.
Shnider and J. Stasheff \cite{MSS} gives a full treatment of
this notion.\\

An \emph{algebraic operad} is an algebraic object that models the
operations acting on certain type of algebras. For instance, there
exists an operad $\A$ coding the operations of associative
algebras, an operad $\C$ for commutative algebras and an operad
$\Li$ for Lie
algebras.\\

The operations acting on $n$ variables form a right $\Sn$-module
$\Po(n)$. A collection $(\Po(n))_{n\in \mathbb{N}^*}$ of
$\Sn$-modules is called an \emph{$\Sy$-module}. One defines a
monoidal product $\circ$ in the category of $\Sy$-modules by the
formula :

$$\Po \circ \Qo (n) := \bigoplus_{1\leqslant k \leqslant n} \left(
\bigoplus_{i_1+\cdots+i_k=n} \Po(k) \otimes \big( \Qo(i_1) \otimes
\cdots \otimes \Qo(i_k)\big) \otimes_{\Sy_{i_1}\times \cdots
\times \Sy_{i_k}} k[\Sy_n]  \right)_{\Sy_k},$$

where the coinvariants are taken with respect to the action of the
symmetric group $\Sy_k$ given by $(p\otimes q_1 \ldots q_k\otimes
\sigma)^\nu := p^\nu \otimes q_{\nu(1)} \ldots q_{\nu(k)} \otimes
\bar{\nu}^{-1}.\sigma$ for $p\in \Po(k)$, $q_j\in \Qo(i_j)$,
$\sigma\in \Sy_n$ and $\nu\in \Sy_k$, such that $\bar{\nu}$ is the
induced block permutation.\\

\begin{figure}[h]
$$\xymatrix@R=20pt@C=20pt{ a_1 \ar[dr] & a_2 \ar[d] & a_3 \ar[dl]& a_4\ar[dr]& &a_5\ar[dl]
& a_6\ar[dr] & a_7\ar[d] & a_8\ar[dl] \\
\ar@{.}[r] & *+[F-,]{q_1}\ar@{.}[rrr] \ar[drrr]& & &
*+[F-,]{q_2}\ar@{.}[rrr] \ar[d]
 & & & *+[F-,]{q_3} \ar@{.}[r]
\ar[dlll]& \\
\ar@{.}[rrrr] & & & & *+[F-,]{p_1} \ar[d]\ar@{.}[rrrr] & & & & \\
 & & & & & & & &  }$$
\caption{\label{fig1} Example of composition of operations
$\Po\circ \mathcal{Q}(8)$}
\end{figure}
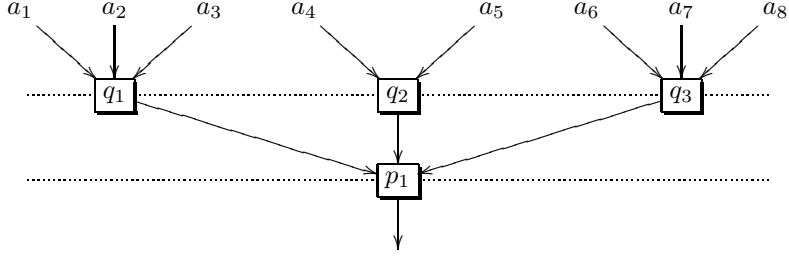

This product reflects the compositions of operations and an
element of $\Po \circ \Qo$ can be represented by 2-levelled trees
whose vertices are indexed by the elements of the operads
(\emph{cf.} [Figure~\ref{fig1}]).

The unit of this monoidal category is given by the $\Sy$-module
$I=(k,\, 0,\, 0,\, \ldots )$.

\begin{dei}[Operad]
An \emph{operad} (or \emph{algebraic operad}) $(\Po,\, \mu,\,
\eta)$ is a monoid in the monoidal category $(\Sy\textrm{-Mod},\,
\circ,\, I)$. This means that the \emph{composition} morphism
$\mu\, : \, \Po \circ \Po \to \Po$ is associative and that the
morphism $\eta \, :\, I \to \Po$ is a \emph{unit}.
\end{dei}

\begin{ex1}
Let $V$ be a $k$-module. The $\Sy$-module $(Hom_k(V^{\otimes n},\
V))_{n\in \mathbb{N}^*}$ of morphisms between $V^{\otimes n}$ and
$V$ forms an operad with the classical composition of morphisms.
This operad is denoted by $End(V)$.
\end{ex1}

We have recalled the definition of an operad in the symmetric
category of $k$-module. One can generalize it to any symmetric
category. For instance, we consider the category $(Set, \times)$
of $Sets$ equipped with the symmetric monoidal product given by
the cartesian product $\times$. We call an \emph{$\Sy$-Set} a
collection $(\Po_n)_{n\in \mathbb{N}^*}$ of sets $\Po_n$ equipped
with an action of the group $\Sy_n$.

\begin{dei}[Set operad]
A monoid $(\Po,\, \mu,\, \eta)$ in the monoidal category of
$\Sy$-sets is called a \emph{set operad}.
\end{dei}

If $\Po$ is a set operad, then the free $k$-module
$\widetilde{\Po}(n):=k[\Po_n]$ is an (algebraic) operad.\\

For any element $(\nu_1, \ldots,\, \nu_t)$ of $\Po_{i_1}\times
\ldots \times \Po_{i_t}$, we denote by $\mu_{(\nu_1, \ldots,\,
\nu_t)}$ the following map defined by the product of the set
operad $\Po$.
\begin{eqnarray*}
\mu_{(\nu_1,\ldots,\, \nu_t)} \ : \ \Po_t &\to &
\Po_{i_1+\cdots+i_t} \\
 \nu  &
\mapsto & \mu(\nu(\nu_1,\ldots ,\, \nu_t)).
\end{eqnarray*}

\begin{dei}[Basic-set operad]
A \emph{basic-set} operad is a set operad $\Po$ such that for any
$(\nu_1,\ldots,\, \nu_t)$ in $(\Po_{i_1},\ldots ,\,\Po_{i_t})$ the
composition maps $\mu_{(\nu_1,\ldots,\, \nu_t)}$ are injective.
\end{dei}

\begin{ex}
The operads $\C$, $\Pe$ and $\A$, for instance, are basic-set
operads (\emph{cf.} section~\ref{ExamplesPosets}).
\end{ex}

This extra condition will be crucial in the proof of the
theorem~\ref{ThmOrdNor}.

\begin{rem}
Dually, one can define the notion of \emph{cooperad} which is an
operad in the opposite category $(\Sy\textrm{-mod}^{\textrm{op}},
\, \circ ,\, I)$. A cooperad $(C,\, \Delta,\, \epsilon)$ is an
$\Sy$-module $C$ equipped with a map $\Delta \, : \, C \to C \circ
C$ coassociative and a map $\epsilon \, : \, C \to I$ which is a
counit.
\end{rem}

An operad $(\Po,\, \mu,\, \eta)$ is \emph{augmented} if there
exists a morphism of operads $\epsilon \, : \, \Po \to I$
such that $\epsilon \circ \eta=id_I$.

\begin{dei}[$\Po$-algebra]
A \emph{$\Po$-algebra} structure on a module $V$ is given by a
morphism $\Po \to End(V) $ of operads.
\end{dei}

It is equivalent to have a morphism $$\mu_A \  : \quad
\Po(A):=\bigoplus_{n\in \mathbb{N}^*} \Po(n)\otimes_{\Sn}
A^{\otimes n} \to A$$ such that the following diagram commutes
$$\xymatrix{\Po\circ\Po (A) \ar[r]^{\mu(A)} \ar[d]^{\Po(\mu_A)} & \Po(A) \ar[d]^{\mu_A}\\
\Po(A) \ar[r]^{\mu_A} & A.} $$

\begin{rem}
The free $\Po$-algebra on the module $V$ is defined by the module
$\Po(V)=\bigoplus_{n\in \mathbb{N}^*} \Po(n) \otimes_{\Sy_n}
V^{\otimes n}$.
\end{rem}

\subsubsection{\bf Free and quadratic operads}

\begin{dei}[Free operad, $\F(V)$]
The forgetful functor from the category of operads to the category
of $\Sy$-modules has a left adjoint functor which gives the
\emph{free operad} on an $\Sy$-module $V$, denoted by $\F(V)$.
\end{dei}

\begin{rem}
The construction of the free operad on $V$ is given by trees whose
vertices are indexed by the elements of $V$. The free operad is
equipped with a natural graduation which corresponds to the number
of vertices of the trees. We denote this graduation by
$\F_{(n)}(V)$.

On the same $\Sy$-module $\F(V)$, one can define two maps $\Delta$
and $\epsilon$ such that $(\F(V),\, \Delta,\, \epsilon)$ is the
cofree connected cooperad on $V$. This cooperad is denoted by
$\F^c(V)$.

For more details on free and cofree operads, we refer the reader
to \cite{Vallette}.
\end{rem}

\begin{dei}[Quadratic operad]
A \emph{quadratic operad $\Po$} is an operad $\Po=\F(V)/(R)$
generated by an $\Sy$-module $V$ and a space of relations
$R\subset \F_{(2)}(V)$, where $\F_{(2)}(V)$ is the direct summand
of $\F(V)$ generated by the trees with $2$ vertices.
\end{dei}

\begin{ex}
$ $
\begin{itemize}
\item[$\triangleright$] The operad $\A$ coding associative
algebras is a quadratic operad generated by a binary operation
$V=\mu.k[\Sy_2]$ and the associative relation $R=\big(
\mu(\mu(\bullet,\, \bullet),\, \bullet)-\mu(\bullet,\,
\mu(\bullet,\, \bullet)\big).k[\Sy_3]$

\item[$\triangleright$] The operad $\C$ coding commutative
algebras is a quadratic operad generated by a symmetric binary
operation $V=\mu.k$ and the associative relation $R=
\mu(\mu(\bullet,\, \bullet),\, \bullet)-\mu(\bullet,\,
\mu(\bullet,\, \bullet)$

\item[$\triangleright$] The operad $\Li$ coding Lie algebras is a
quadratic operad generated by an anti-symmetric operation
$V=\mu.sgn_{\Sy_2}$ and the Jacobi relation $R=\mu(\mu(\bullet,\,
\bullet),\, \bullet). (id+(123)+(132))$
\end{itemize}
\end{ex}

Throughout the text, we will only consider quadratic operads generated by
$\Sy$-modules $V$ such that $V(1)=0$. Therefore, the operad verifies
$\Po(1)=k$. We say that $V$ is an \emph{homogenous} $\Sy$-module
if $V(t)\ne 0$ for only one $t$.

\subsection{Koszul Duality for operads}
We recall the results of the theory of Koszul duality for operads
that will be used in this article. This theory was settled by V.
Ginzburg and M.M. Kapranov in \cite{GK}. The reader can find a
short survey of the subject in the article of J.-L. Loday
\cite{LodayBourbaki} and in the appendix B of \cite{Loday}. For a
full treatment and the generalization over a field of
characteristic $p$ (or a Dedekind ring), we refer to the article
of B. Fresse \cite{Fresse}.

\subsubsection{\bf Koszul dual operad and cooperad}

To a quadratic operad, one can associate a Koszul dual operad and
cooperad.\\

The Czech dual of an $\Sy$-module $V$ is the linear dual of $V$
twisted by the signature representation.

\begin{dei}[Czech dual of an $\Sy$-module, $V^\vee$]
To an $\Sy$-module $V$, one can associate another $\Sy$-module,
called the \emph{Czech dual of $V$} and denoted by $V^\vee$, by
the formula $V^\vee(n)=V(n)^*\otimes sgn_{\Sn}$, where $V^*$ is
the linear dual of $V$ and $sgn_{\Sn}$ is the signature
representation of $\Sn$.
\end{dei}

\begin{dei}[Koszul dual of an operad, $\Po^!$]
Let $\Po=\F(V)/(R)$ be a quadratic operad. Its \emph{Koszul dual
operad} is the quadratic operad generated by the $\Sy$-module
$V^\vee$ and the relation $R^\perp$. This operad is denoted by
$\Po^!=\F(V^\vee)/(R^\perp)$.
\end{dei}

\begin{ex}
$ $
\begin{itemize}
\item[$\triangleright$] The operad $\A$ is autodual $\A^!\cong\A$.

\item[$\triangleright$] The operads $\C$ and $\Li$ are dual to
each other, $\C^!\cong \Li$ and $\Li^!\cong \C$.

\end{itemize}
\end{ex}

\begin{dei}[Koszul dual cooperad, $\Po^{\ac}$]
Let $\Po=\F(V)/(R)$ be a quadratic operad generated by a finite
dimensional $\Sy$-module $V$. The \emph{Koszul dual cooperad} of
$\Po$ is defined by the Czech dual of $\Po^!$ and denoted by
$\Po^{\ac}={\Po^!}^\vee$.
\end{dei}

\subsubsection{\bf Koszul complex and Koszul operads}
An operad is said to be a \emph{Koszul} operad if its Koszul
complex is acyclic.

\begin{dei}[Koszul complex]
To a quadratic operad $\Po$, one associates the following boundary
map on the $\Sy$-module $\Po^{\ac}\circ \Po$ :

$$\partial \ : \  \Po^{\ac}\circ \Po \xrightarrow{\Delta'} \Po^{\ac}\circ \underbrace{\Po^!}_1 \circ
\, \Po \xrightarrow{\Po^{\ac} \circ (I+\alpha) \circ \Po}
\Po^{\ac}\circ\Po \circ \Po \xrightarrow{\Po^{\ac}\circ \mu_{\Po}}
\Po^{\ac}\circ \Po,$$ where $\Delta'$ is the projection of the
coproduct $\Delta={\mu_{\Po^!}}^\vee$ of the cooperad $\Po^{\ac}$
on the $\Sy$-module $\Po^!\circ \underbrace{\Po^!}_1$ with only
one element of $\Po^!$ on the right and $\alpha$ is the following
composition :
$$\alpha : \Po^{\ac}=(\F(V^\vee)/(R^\perp))^\vee \twoheadrightarrow
(V^\vee)^\vee \simeq V \hookrightarrow \Po=\F(V)/(R).$$
\end{dei}

\begin{dei}[Koszul operad]
A quadratic operad $\Po$ is a \emph{Koszul} operad if its Koszul
complex $(\Po^{\ac}\circ \Po,\, \partial)$ is acyclic.
\end{dei}

\begin{ex}
The operads $\A$, $\C$ and $\Li$ are Koszul operads.
\end{ex}

\subsubsection{\bf Operadic homology}
To any quadratic operad $\Po=\F(V)/(R)$, one can define a homology
theory for $\Po$-algebras.

\begin{dei}[Operadic chain complex, $\Co^\Po_*(A)$]
Let $A$ be $\Po$-algebra. On the modules
$\Co^\Po_n(A):=\Po^{\ac}(n)\otimes_{\Sn}A^{\otimes n}$, one
defines a boundary map by the formula :
$$\Po^{\ac}(A) \xrightarrow{\Delta'} \Po^{\ac}\circ \underbrace{\Po^!}_1 \circ A
\xrightarrow{\Po^{\ac} \circ (I+\alpha) \circ A}
\Po^{\ac}\circ\Po(A) \xrightarrow{\Po^{\ac}\circ\mu_A}
\Po^{\ac}(A).$$
The related homology is denoted by $H^\Po_*(A)$.
\end{dei}

\begin{ex}
$ $
\begin{itemize}

\item[$\triangleright$] In the $\A$ case, the operadic homology
theory correspond to the Hochschild homology of associative
algebras.

 \item[$\triangleright$] In the $\C$ case, the operadic homology
theory correspond to the Harrison homology of commutative
algebras.

\item[$\triangleright$] In the $\Li$ case, the operadic homology
theory correspond to the Chevalley-Eilenberg homology of Lie
algebras.
\end{itemize}
\end{ex}

Let us recall that an operad $\Po$ is a Koszul operad if the
operadic homology of the free $\Po$-algebra $\Po(V)$ is acyclic
for every $k$-module $V$,
$$H^\Po_*(\Po(V))=\left\{
\begin{array}{l}
V \quad \textrm{if} \quad *=1,
\\ 0 \quad \textrm{elsewhere}.
\end{array} \right.$$

\begin{rem}
An operad $\Po$ is Koszul if and only if its dual $\Po^!$ is
Koszul.
\end{rem}

\subsubsection{\bf Differential bar construction of an operad}

One can generalize the differential bar construction of
associative algebras to operads. \\

Let $\Po=\F(V)/(R)$ be a quadratic operad. The natural graduation
of the free operad $\F(V)$ induces a graduation on $\Po$. We call
it the \emph{weight} and we denote it by $\Po_{(n)}:=\textrm{Im}
\left( \F_{(n)}(V)\to \F(V)/(R) \right)$. Any quadratic operad is
augmented and the ideal of augmentation is given
by $\bar{\Po}:=\bigoplus_{n > 0} \Po_{(n)}$.\\

We consider the \emph{desuspension} of the ideal of augmentation $\Sigma^{-1}
\bar{\Po}$ which corresponds to shift the weight by $-1$. An element $\nu$
of weight $n$ in $\bar{\Po}$ gives an element $\Sigma^{-1}\nu$ of weight
$n-1$ in $\Sigma^{-1}\bar{\Po}$.

\begin{dei}[Differential bar construction of an operad, $\B(\Po)$]
The \emph{differential bar construction} of a quadratic operad $\Po$
is given by the cofree connected cooperad on $\Sigma^{-1} \bar{\Po}$, namely
$\F^c(\Sigma^{-1} \bar{\Po})$, where the coboundary map
is the unique coderivation
$\delta$ induced by the partial product of $\Po$ :
$$\F^c_{(2)}(\Sigma^{-1} \bar{\Po}) \xrightarrow{\bar{\mu}}
\Sigma^{-1} \bar{\Po}.$$

This differential
cooperad is denoted by $\left(\B(\Po), \delta \right)$.
\end{dei}

The cohomological degree of the differential bar construction is
induced by the global weight of $\F^c(\Sigma^{-1} \bar{\Po})$ and
the desuspension. We denote it by $\B^n(\Po)$. For instance, a
tree with $l$ vertices indexed by operations $\Sigma^{-1}\nu_1,
\ldots,\, \Sigma^{-1}\nu_l$ of weight $p_1-1,\ldots,\, p_l-1$
represents an element of cohomological degree equals to
$p_1+\cdots+p_l-l$.

\subsubsection{\bf Koszul duality}

We recall here the main theorem of the Koszul duality theory for operads. \\

For a quadratic operad, the differential bar construction has the following form :

$$ 0 \xrightarrow{\delta}
\B^{0}(\Po)=\F^c(\Sigma^{-1}V)
  \xrightarrow{\delta}
\B^{1}(\Po)
 \xrightarrow{\delta}  \cdots $$

\begin{rem}
The differential bar construction is a
direct sum of subcomplexes indexed by the global weight coming from
$\Po$ taken without the desuspension. We denote it by
$\B(\Po)=\bigoplus_{m\in \mathbb{N}} \B(\Po)_{(m)}$.
\end{rem}

The main theorem of the Koszul duality theory used in this
article is the following one.

\begin{thm}{\rm [GK-F]}
\label{GK-F}
Let $\Po$ be a quadratic operad generated by a finite dimensional space $V$.
\begin{enumerate}
\item One has $H^0(\B(\Po))\cong \Po^{\ac}$ (more precisely, one has
$H^0(\B(\Po)_{(m)})\cong \Po^{\ac}_{(m)}=(\Po^!_{(m)})^*)$.
\item The operad $\Po$ is a Koszul operad if and only if $H^k(\B^*(\Po))=0$ for
$k>0$.
\end{enumerate}
\end{thm}

\begin{rem}
Notice the similarity with the notion of a Cohen-Macaulay poset. A
Cohen-Macaulay poset has its homology concentrated in top
dimension and a Koszul operad is an operad such that the homology
of its bar construction is concentrated in some degree too.
\end{rem}

\section{Partition posets associated to an operad}

We give the definition of the partition poset associated to an
operad. We recall the definition of the simplicial bar
construction of an operad. And we prove that the operad is Koszul
if and only if the related poset is Cohen-Macaulay.

\subsection{Construction of a partition poset from an operad}
To any set operad, we define a family of posets of partitions.

\subsubsection{\bf $\Po$-partition}

Let $(\Po,\, \mu,\, \eta)$ be a set operad. By definition, the set
$\Po_n$ is equipped with an action of the symmetric group $\Sy_n$.
Let $I$ be a set $\{x_1,\ldots,\, x_n\}$ of $n$ elements. We
consider the cartesian product $\Po_n \times \mathcal{I}$, where
$\mathcal{I}$ is the set of ordered sequences of elements of $I$,
each element appearing once. We define the diagonal action of
$\Sy_n$ on $\Po_n\times \mathcal{I}$ as follows. The image of
$\nu_n\times (x_{i_1}, \ldots,\, x_{i_n})$ under a permutation
$\sigma$ is given by $\nu_n.\sigma \times
(x_{\sigma^{-1}(i_1)},\ldots,\, x_{\sigma^{-1}(i_n)})$. We denote
by $\overline{\nu_n\times(x_{i_1},\ldots ,\, x_{i_n})}$ the orbit
of $\nu_n\times(x_{i_1},\ldots ,\, x_{i_n})$ and by
$\Po_n(I):=\Po_n\times_{\Sy_n} \mathcal{I}$ the set of orbits
under this action. We have
$$\Po(I)=\left( \coprod_{f\, : \, \textrm{bijection} \atop [n] \to I}
\Po_n\right)_\sim,$$ where the equivalence relation $\sim$ is
given by $(\nu_n, \, f)\sim (\nu_n. \sigma, \, f\circ
\sigma^{-1})$.

\begin{dei}[$\Po$-partition]
Let $\Po$ be a set operad. A \emph{$\Po$-partition} of $[n]$ is a
set of components $\{B_1,\ldots,\, B_t \}$ such that each $B_j$
belongs to $\Po_{i_j}(I_j)$ where $i_1+\cdots+i_t=n$ and
$\{I_j\}_{1 \leqslant j \leqslant t}$ is a partition of $[n]$.
\end{dei}

\begin{rem}
A $\Po$-partition is a classical partition enriched by the
operations of $\Po$.
\end{rem}

\subsubsection{\bf Operadic partition poset} We give the
definition of the poset associated to an operad.\\

We generalize the maps $\mu_{(\nu_1,\ldots,\,\nu_t)}$ (see section
$2$) to the $\Po(I)$.

\begin{pro}
Let $\{ I_1, \ldots,\, I_t \}$ be a partition of a set $I$, where
the number of elements of $I_r$ is equal to $i_r$. Let $(C_1,
\ldots, C_t)$ be an element of $\Po_{i_1}(I_1)\times \cdots \times
\Po_{i_t}(I_t)$. Each $C_r$ can be represented by
$\overline{\nu_r\times(x^r_1,\ldots ,\, x^r_{i_r})}$, where
$\nu_r\in\Po_{i_r}$ and $I_r=\{ x^r_1,\ldots ,\, x^r_{i_r} \}$.
The map $\tilde{\mu}$ given by the formula

\begin{eqnarray*}
\tilde{\mu} \ : \  \Po_t\times \left( \Po_{i_1}(I_1)\times \cdots
\times  \Po_{i_t}(I_t) \right)
 &\to &
\Po_{i_1+\cdots+i_t}(I) \\
\nu\times \left( C_1
, \ldots , C_r
 \right) &
\mapsto &
\overline{\mu(\nu(\nu_1,\ldots ,\, \nu_t))\times (x^1_1,\ldots ,\, x^t_{i_t})},
\end{eqnarray*}
is well defined and equivariant under the action of $\Sy_t$.
\end{pro}

\begin{deo}
It is a direct consequence of the definition of a set operad.
$\cqfd$
\end{deo}

To a set operad $\Po$, we associate a partial order on the
set of $\Po$-partitions as follows.

\begin{dei}[Operadic partition poset, $\Pi_\Po$]
Let $\Po$ be a set operad.

Let $\lambda=\{B_1,\ldots ,\, B_r\}$ and $\mu=\{C_1,\ldots ,\, C_s
\}$ be two $\Po$-partitions of $[n]$, where $B_k$ belongs to
$\Po_{i_k}(I_k)$ and $C_l$ to $\Po_{j_l}(J_l)$. The
$\Po$-partition $\mu$ is said to be \emph{larger} than $\lambda$
if, for any $k\in \{1,\ldots,\, r\}$, there exist $\{p_1,
\ldots,\, p_t\} \subset \{1,\ldots ,\, s\}$ such that
$\{J_{p_1},\ldots ,\, J_{p_t}\}$ is a partition of $I_k$ and if
there exists an element $\nu$ in $\Po_t$ such that
$B_k=\tilde{\mu}(\nu\times (C_{p_1},\ldots,\, C_{p_t}))$. We
denote this relation by $\lambda \leqslant \mu$.

We call this poset, the \emph{operadic partition poset} associated to the operad
$\Po$ and we denote it by $\Pi_{\Po}$.
\end{dei}

Throughout the text, we suppose that the operad $\Po$ is such that
the set $\Po_1$ is reduced to one element, the identity for the
composition. In this context, the poset $\Pi_{\Po}(n)$ has only
one maximal element corresponding to the partition $\{ \{ 1\},
\ldots ,\, \{ n\}\}$, where $\{ i\}$ represents the unique element
of $\Po_1(\{ i\})$. Following the classical notations, we denote
this element by $\hat{1}$. The set of minimal elements is
$\Po_n([n])$.

\begin{rem}
If the operad $\Po$ is a quadratic operad generated by an
homogenous $\Sy$-set $V$ such that $V(t)\neq 0$, then the
$\Po$-partitions have restricted block size. The possible lengths
for the blocks are $l(t-1)+1$ with $l\in \mathbb{N}$.
\end{rem}

\subsubsection{\bf Operadic order complex}

We consider the following order complex associated to the poset
$\Pi_\Po$.

\begin{dei}[Operadic order complex, $\Delta(\Pi_\Po)$]
The \emph{operadic order complex} is the presimplicial complex
induced by the chains $\lambda_0 < \lambda_1 < \cdots < \lambda_l$
of $\Pi_\Po$, where $\lambda_0$ is a minimal element and
$\lambda_l=\hat{1}$. We denoted by $\Delta_l(\Pi_\Po)$ the set of
chains of length $l$.
\end{dei}

\begin{pro}
\label{grsubpos} Let $\Po$ be a set operad.
\begin{enumerate}
\item Each $k$-module $k[\Po_n]$ has finite rank if and only if
each poset $\Pi_\Po(n)$ is finite.

\item Moreover, if $\Po$ is a quadratic operad generated by an homogenous
$\Sy$-set $V$ such that $V(t)\ne 0$. then all the maximal chains of $\Pi_\Po$ have the
same length. More precisely, if $\lambda\leqslant \mu$ the maximal
chains in the closed interval $[\lambda,\, \mu]$ have length
$(b_\mu-b_\lambda)/(t-1)$, where $b_\lambda$ denotes the number of
blocks of the $\Po$-partition $\lambda$.

\noindent The subposets $[\alpha, \, \hat{1}]$, for $\alpha\in
Min(\Pi_\Po(n))=\Po_n([n])$ are graded posets. The length of
maximal chains between a minimal element and $\hat{1}$ in
$\Pi_\Po(n)$ is equal to $l+1$ if $n=l(t-1)+1$ and $0$ otherwise.
\end{enumerate}
\end{pro}

\begin{deo}
The first point is obvious. In the second case, the condition on
$V$ implies that all the blocks have size $l(t-1)+1$. The
surjection of the compositions gives that every block of size
$l(t-1)+1$ is refinable if $l>1$. Since each closed interval of
the form $[\alpha, \, \hat{1}]$, for $\alpha\in
Min(\Pi_\Po(n))=\Po_n([n])$ is bounded and pure, it is graded by
definition. $\cqfd$
\end{deo}

\subsection{Simplicial and normalized bar construction over a Koszul operad}

We recall the definition of the simplicial and the normalized bar
construction of an operad. They are always quasi-isomorphic. We
recall the quasi-isomorphism between the simplicial and the
differential bar construction of an operad due to B. Fresse in
\cite{Fresse}.

\subsubsection{\bf Definition of the simplicial and the normalized
bar construction of an operad}

Following the classical methods of monoidal categories (\emph{cf.}
S. Mac Lane \cite{MacLane}), we recall the construction of the
simplicial bar construction of an operad.

\begin{dei}[Simplicial bar construction of an operad, $\Co(\Po)$]
Let $(\Po,\, \mu,\, \eta,\, \epsilon)$ be an augmented operad. We define the
$\Sy$-module $\Co(\Po)$ by the formula $\Co_l(\Po)=\Po^{\circ l}$.

\begin{tabbing}
The face maps $d_i$ are given by \\
$\quad$ \= $d_i=\Po^{\circ(i-1)}\circ\mu\circ\Po^{\circ (l-i-1)} \ : $ \\
\> $\Po^{\circ l}=\Po^{\circ(i-1)}\circ (\Po \circ \Po)
\circ\Po^{\circ (l-i-1)}\to \Po^{\circ(i-1)}\circ\Po
\circ\Po^{\circ (l-i-1)}=\Po^{\circ (l-1)},$ \\
for $0 < i <l$ and otherwise by \\
\> $d_0=\epsilon \circ \Po^{\circ (l-1)} \ : \ \Po^{\circ l}
\to \Po^{\circ (l-1)}$ \\
\> $d_l=\Po^{\circ (l-1)} \circ \epsilon \ : \ \Po^{\circ l}
\to \Po^{\circ (l-1)}.$ \\
The degeneracy maps  $s_j$ are given by \\
\> $s_j=\Po^{\circ j}\circ\eta\circ\Po^{\circ (l-j)} \ : $ \\
\> $\Po^{\circ l}=\Po^{\circ j}\circ I \circ\Po^{\circ (l-j)}\to
\Po^{\circ j}\circ\Po \circ\Po^{\circ (l-j)}=\Po^{\circ (l+1)},$\\
for $0\leq j \leq l$.
\end{tabbing}

This simplicial module is called the \emph{simplicial (or
categorical) bar construction of the operad $\Po$}.
\end{dei}

Since the composition and unit morphisms of an operad preserve the
action of the symmetric groups $\Sn$, the face and the degeneracy
maps of the simplicial bar construction $\Co(P)$ are morphisms
of $\Sy$-modules.\\

The $\Sn$-modules $\Co_l(\Po)(n)$ can be represented by
$l$-levelled trees with $n$ leaves, whose vertices are indexed by
operations of $\Po$. In this framework, the faces $d_i$ correspond
to the notion of \emph{level contraction}.\\

We consider the normalization of this simplicial complex.

\begin{dei}[Normalized bar construction, $\mathcal{N}(\Po)$]
The \emph{normalized bar construction} of an augmented operad
$\Po$ is the quotient
$$\mathcal{N}_l(\Po):=\mathcal{C}_l(\Po)
\big/ {\textstyle   \sum_{0\leqslant j \leqslant l-1}} \textrm{Im}(s_j)$$
equipped
with the boundary map induced by the faces.
\end{dei}

The $\Sy$-module $\N_l(\Po)$ is isomorphic to the direct summand
of $\Co_l(\Po)$ composed by $l$-levelled trees with at least one
element of $\bar{\Po}$ on each level. The canonical projection of
complexes $\mathcal{C}(\Po) \twoheadrightarrow \mathcal{N}(\Po)$
is a quasi-isomorphism (see
 E.B. Curtis \cite{Curtis} Section $3$ or J.P. May \cite{May}
Chapter V for this classical result).\\

Once again, the simplicial bar construction of a quadratic operad
$\Po$ is a direct sum, indexed by $m\in \mathbb{N}$, of the
subcomplexes of $\mathcal{C}(\Po)$ composed by elements of global
weight $m$, and denoted by $\mathcal{C}(\Po)_{(m)}$. The same
decomposition holds for the normalized bar construction
$\mathcal{N}(\Po)= \bigoplus_{m \in \mathbb{N}}
\mathcal{N}(\Po)_{(m)},$ and the canonical projection preserves
this weight. Therefore, the canonical projection maps
$\mathcal{C}(\Po)_{(m)} \twoheadrightarrow \mathcal{N}(\Po)_{(m)}$
are quasi-isomorphisms.

\subsubsection{\bf Homology of the normalized bar construction and
Koszul duality of an operad} We recall the quasi-isomorphism
between the differential and the normalized bar constructions of
an operad due to B. Fresse. It shows that a quadratic operad is
Koszul if and only if the homology of the normalized bar
construction is concentrated in top dimension. \\

One can consider the homological degree of the bar construction
given by the number of vertices of the trees. We denote it by
$\B_{(l)}(\Po)=\F^c_{(l)}(\Sigma^{-1}\bar{\Po})$.
It is equal to the number of signs $\Sigma^{-1}$ used to write
an element. Be careful, this degree is
different from the cohomological degree $(p_1+\cdots+p_l-l)$ and
also different from the global weight taken without the desuspension
$(p_1+\cdots+p_l)$. The relation between these three graduations is obvious.
The sum of the homological degree with the cohomological degree gives the global
weight. \\

In \cite{Fresse} (section $4$), B. Fresse defines a morphism,
called the \emph{levelization morphism}, between the differential
bar construction and the normalized bar construction of an operad.
It is defined from $\B_{(l)}(\Po)$ to $\mathcal{N}_{l}(\Po)$. This
levelization morphism induces a morphism of complexes when one
considers the homological degree on the differential bar
construction.

\begin{thm}{\rm [F] (Theorem 4.1.8.)}
\label{thmF} \label{qiwg} Let $\Po$ be an operad such that
$\Po(0)=0$ and $\Po(1)=k$. The levelization morphism $\B(\Po) \to
\mathcal{N}(\Po)$ is a quasi-isomorphism.
\end{thm}

\begin{rem}
Since the normalized and the simplicial bar constructions are
quasi-isomorphic, the bar construction $\B(\Po)$ of an operad is
quasi-isomorphic to the simplicial bar construction $\Co(\Po)$.
\end{rem}

As we have seen before, the normalized bar construction of a
quadratic operad $\Po$ is a direct sum, indexed by the weight
$(m)$, of the subcomplexes of $\mathcal{N}(\Po)$, composed by
elements of global weight $m$. The levelization morphism preserves
this weight. Therefore, the legalization morphism
$\B_{(l)}(\Po)_{(m)} \to \mathcal{N}_l(\Po)_{(m)}$ is a
quasi-isomorphism for any $l$ and any $m$.

\begin{cor}
\label{coro}
Let $\Po$ be a quadratic operad generated by a finite dimensional
space $V$.
\begin{enumerate}
\item We have $H_m(\mathcal{N}_*(\Po)_{(m)})\cong \Po^{\ac}_{(m)}.$

\item The operad $\Po$ is a Koszul operad if and only if
$H_l(\mathcal{N}_*(\Po)_{(m)})=0$ for $l \ne m$.
\end{enumerate}
\end{cor}

\begin{deo}
This corollary is the analog of Theorem~\ref{GK-F} with the
quasi-isomorphisms given in Theorem~\ref{qiwg}. For the first
point, we have the following identities
$$\Po^{\ac}_{(m)}\cong H^0(\B^*(\Po)_{(m)}) = H_m(\B_{(*)}(\Po)_{(m)})
\cong H_m(\mathcal{N}_{*}(\Po)_{(m)}).$$ For the second point, we
know from Theorem~\ref{GK-F} that $\Po$ is a
 Koszul operad if and only if
$H^k(\B^*(\Po))=0$ for any $k>0$, which is equivalent to
$H_l(\B_{(*)}(\Po){(m)})=0$ for $l\ne m$ with the homological
degree. The isomorphisms between the differential and the
normalized bar constructions of Theorem~\ref{qiwg} give that $\Po$
is a Koszul operad is and only if
$H_l(\mathcal{N}_{(*)}(\Po)_{(m)})=0$ for $l\ne m$. $\cqfd$
\end{deo}

\subsection{Homology of a partition poset associated to an operad}
We prove that an operad is Koszul if and only if the related
complex is Cohen-Macaulay. In this case, the homology of the poset
is isomorphic, as an $\Sy$-module, to the Koszul dual cooperad.

Let $\Po$ be a set operad and $\widetilde{\Po}=k[\Po]$ the corresponding algebraic operad.

\begin{figure}[h]
$$\xymatrix@R=20pt@C=20pt{ &  {4} \ar@{-}[dr] & {2} \ar@{-}[d] & {5} \ar@{-}[d] &
{1}\ar@{-}[d] &
  {7} \ar@{-}[dr] & {6} \ar@{-}[d] & {3} \ar@{-}[dl]  \\
 *+[o][F-]{3} &  \ar@{.}[r]  &*+[F-,]{\nu_1} \ar@{-}[dr]
&*+[F-,]{\nu_2}\ar@{.}[l]\ar@{-}[d] &*+[F-,]{\nu_3}
\ar@{-}[dl]\ar@{.}[l] \ar@{.}[rr]&  &   *+[F-,]{\nu_4} \ar@{-}[dl]
\ar@{.}[r]
 & \\
*+[o][F-]{2} &\ar@{.}[rr] &  &*+[F-,]{\nu_5} \ar@{-}[dr]\ar@{.}[rr]
& & *+[F-,]{\nu_6}\ar@{-}[dl]\ar@{.}[rr]  & & \\
*+[o][F-]{1} &\ar@{.}[rrr]  &  &  & *+[F-,]{\nu_7} \ar@{-}[d]
\ar@{.}[rrr]  & & & \\
& & & & & & }$$
 \caption{\label{figC} An example of a $3$-levelled tree
 with the vertices indexed by elements of
$\Po$.}
\end{figure}
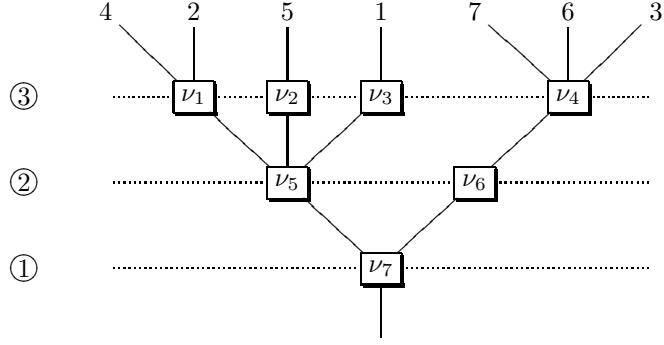

\begin{dei}[$\mathcal{C}(\Po)(I)$ and $\mathcal{N}(\Po)(I)$]
Let $I$ be a set with $n$ element ($\{x_1,\ldots ,\, x_n \}$ or
$[n]$ for instance). The set $\mathcal{C}_l( \Po)(I)$ is the set
of $l$-levelled trees with the vertices indexed by elements of
$\Po$ and where each leaf is indexed by a different element of $I$
(\emph{cf.} figure~\ref{figC}). When $I=[n]$, we denote the set
$\mathcal{C}(\Po)(I)$ by $\mathcal{C}(\Po)(n)$.

The normalized part $\mathcal{N}(\Po)(I)$ is given by the subset
of $\mathcal{C}(\Po)(I)$ generated by the same levelled trees such
that at least one vertex on each level is indexed by an element of
$\Po_n$ with $n>1$.
\end{dei}

\begin{lem}
$ $
\begin{enumerate}
\item The set $\mathcal{C}(\Po)(n)$ is stable under the face maps
$d_i$ and the degeneracy maps $s_j$ of the simplicial bar
construction $\mathcal{C}(\widetilde{\Po})(n)$, making it a simplicial set. The
set $\mathcal{N}(\Po)(n)$ is its normalized associated set.

\item The set $\mathcal{C}(\Po)$ is a basis over $k$ of the simplicial
bar construction $\mathcal{C}(\widetilde{\Po})$ and $\mathcal{N}(\Po)$ is a basis
of the normalized bar construction $\mathcal{N}(\widetilde{\Po})$.
\end{enumerate}
\end{lem}

\begin{deo}
The proof is a direct consequence of the definition of a set operad $\Po$.
$\cqfd$
\end{deo}

In the rest of the text, we will need that the set operad $\Po$ is
a basic set operad (see section $2$).

\begin{thm}
\label{ThmOrdNor} Let $\Po$ be a basic-set operad.

For any $n\in \mathbb{N}^*$, the order set $\Delta_*(\Pi_\Po(n))$ is isomorphic to
$\mathcal{N}_*(\Po)(n)$. This isomorphism preserves the face
maps $d_i$ for $0<i<l$ and the action of  $\Sy_n$. It induces an
isomorphism of presimplicial $k[\Sy_n]$-modules between the order
complex $k[\Delta(\Pi_\Po)(n)]$ and the normalized bar
construction $\mathcal{N}(\widetilde{\Po})(n)$.
\end{thm}

\begin{deo}
We are going to describe a bijection $\Psi$ between
$\mathcal{N}_l(\Po)(n)$
and $\Delta_l(\Pi_\Po(n))$. \\

Let $\mathfrak{T}$ be a non planar tree with $l$ levels and $n$
leaves whose vertices are indexed by elements of $\Po  $. To such
a tree, we build a chain of $\Po$-partitions of $[n]$ in the
following way : we cut the tree $\mathfrak{T}$ along the $i^{th}$
level and we look upwards. We get $t$ indexed and labelled
subtrees. Each of them induces an element in a $\Po  _{i_j}(I_j)$
by composing the operations indexing the vertices along the scheme
given by the subtree. The figure~\ref{FigDemoPsi} shows an example
in the case of the operad $\C$.

\begin{figure}[h]
$$\xymatrix@R=15pt@C=15pt{
  {1} \ar@{-}[dr] & {4} \ar@{-}[d] & {3} \ar@{-}[d] & {6}\ar@{-}[d] &
  {5} \ar@{-}[dr] & {2} \ar@{-}[d] & {7} \ar@{-}[dl] & \\
 *+[o][F-]{3} $ $    \ar@{.}[r]  &*{} \ar@{-}[dr]
&*{}\ar@{.}[l]\ar@{-}[d] &*{}
\ar@{-}[dl]\ar@{.}[l] \ar@{.}[rr]&  &   *{} \ar@{-}[dl] \ar@{.}[r]
 & &  \ar@{<->}[l]{ \lambda_3= \{\{ 1\},\, \{ 4\},\, \{ 3\},\, 
\{ 6\},\, \{ 5\},\, \{ 2\},\,
\{ 7\}  \}}\\
*+[o][F-]{2} \ar@{.}[rr] &  &*{} \ar@{-}[dr]\ar@{.}[rr] &
& *{}\ar@{-}[dl]\ar@{.}[rr]  & &  &
 \ar@{<->}[l]{\lambda_2=\{
\{1,\,4 \},\,  \{ 3\},\, \{ 6\},\, \{2,\,5,\, 7\}
\}}\\
*+[o][F-]{1} \ar@{.}[rrr]  &  &  & *{} \ar@{-}[d]
\ar@{.}[rrr]  & & & &  \ar@{<->}[l]{\lambda_1=
\{ \{1,\,3 ,\,4 ,\,6 \},\, \{ 2,\,5 ,\,7\} \}
}\\
& & & & & & & \ar@{<->}[l] {\lambda_0=
\{ \{1 ,\,2 ,\,3 ,\,4 ,\,5,\, 6,\,7 \}  \}
}}$$
\caption{\label{FigDemoPsi}
An example of the image of $\Psi$ in the case $\Po=\C$.}
\end{figure}

The union of these elements forms a $\Po$-partition
$\lambda_i$ of $[n]$. At the end, by cutting the tree
$\mathfrak{T}$ at the root, we get one minimal $\Po$-partition
$\lambda_0\in \Po_n([n])$.\\

Since $\lambda_{i+1}$ is a strict refinement of $\lambda_i$, the
growing sequence of $\Po$-partitions $\lambda_0 < \lambda_1 <
\cdots < \lambda_l$ is an element of $\Delta_l(\Pi_P(n))$. We
denote this map by $\Psi$. \\

The surjectivity of the map $\Psi$ comes from the definition of
the partial order between the $\Po$-partitions. And the
injectivity of the maps $\mu_{(\nu_1,\ldots,\, \nu_t)}$ induces
the injectivity of $\Psi$. Therefore, $\Psi$ is a bijection.\\

Contracting the $i^{th}$ and the $(i+1)^{th}$ levels of the tree
$\mathfrak{T}$ corresponds, via $\Psi$, to removing the $i^{th}$
partition of the chain $\lambda_0 < \lambda_1 < \cdots <
\lambda_l$. Moreover $\Psi$ preserves the action of the symmetric
group $\Sn$. Therefore, $\Psi$ induces an isomorphism of
presimplicial $k[\Sn]$-modules. $\cqfd$
\end{deo}

\begin{thm}
\label{KoszulPoset} Let $\Po$ be a basic-set and quadratic operad generated by
an homogenous $\Sy$-set $V$ such that $V(t)\neq 0$.

\begin{enumerate}
\item We have $H_m(\Pi_\Po(m(t-1)+1))\cong {\widetilde{\Po}}^{\ac}_{(m)}$.

\item The operad $\widetilde{\Po}$ is a Koszul operad if and only if
$H_l(\Pi_\Po(m(t-1)+1))=0$ for $l\neq m$.
\end{enumerate}
\end{thm}

\begin{deo}
Since the operad $\Po$ is quadratic and generated by $t$-ary
operations, we have
$\mathcal{N}_*(\Po)_{(m)}=\mathcal{N}_*(\Po)(m(t-1)+1)$. The
theorem~\ref{ThmOrdNor} implies that
$H_l(\mathcal{N}_*(\Po)(m(t-1)+1))=H_l(\Pi_\Po(m(t-1)+1)$ for any
$k$. With this quasi-isomorphism, this theorem is a rewriting of
the corollary~\ref{coro}, in the framework of posets.$\cqfd$
\end{deo}

With the notion of Cohen-Macaulay posets, this theorem implies the
following one.

\begin{thm}
\label{KoszulCohen} Let $\Po$ be a basic-set and quadratic operad generated by
an homogenous $\Sy$-set $V$ such that $V(t)\neq 0$.

The operad $\widetilde{\Po}$ is a Koszul operad if and only if
each subposet $[\alpha,\, \hat{1}]$ of each $\Pi_\Po(n)$ is
Cohen-Macaulay, where $\alpha$ belongs to
$Min(\Pi_\Po(n))=\Po_n([n])$.
\end{thm}

\begin{deo}
$ $

$(\Rightarrow)$ The proposition~\ref{grsubpos} shows that
$\Pi_\Po$ is pure. Therefore, each subposet of the form
$[\alpha,\, \hat{1}]$ is graded.

If the operad $\widetilde{\Po}$ is Koszul, we have by
Theorem~\ref{KoszulPoset} that the homology of each poset
$\Pi_\Po(m(t-1)+1)$ is concentrated in top dimension $m$. Since
$$H_l(\Pi_\Po(m(t-1)+1)=\bigoplus_{\alpha \in
Min(\Pi_\Po(m(t-1)+1)} \widetilde{H}_{l-2}(\,(\alpha,\, \hat{1})\,
),
$$ we have $H_l([\alpha,\, \hat{1}]) = \widetilde{H}_{l-2}(\, (\alpha,\, \hat{1}) \,)=0$
for $l\neq m$. Let $x\leqslant y$ be two elements of
$\Pi_\Po(m(t-1)+1)$. Denote the $\Po$-partition $x$ by $\{
B_1,\ldots ,\, B_r\}$ and $y$ by $\{ C_1,\ldots ,\, C_s\}$. Each
$B_k\in \beta_{i_k}(I_k)$ is refined by some $C_l$. For any
$1\leqslant k \leqslant r$, we consider the subposet $[x_k,\,
y_k]$ of $\Pi_\Po(I_k)$, where $x_k=B_k$ and $y_k$ the
corresponding set of $C_l$. There exists one $\alpha_k \in
Min(\Pi_\Po(|y_k|))$ such that the poset $[x_k,\, y_k]$ is
isomorphic to $[\alpha_k,\, \hat{1}]$, which is a subposet of
$\Pi_\Po(|y_k|)$. (The notation $|y_k|$ represents the number of
$C_l$ in $y_k$.) This decomposition gives, with the K\"unneth
theorem, the following formula

\begin{eqnarray*}
\widetilde{H}_{l-1}(\, (x,\, y) \, ) \cong
\bigoplus_{l_1+\cdots+l_r=l} \bigotimes_{k=1}^r
\widetilde{H}_{l_k-1} (\, (x_k,\, y_k) \, ) \cong
\bigoplus_{l_1+\cdots+l_r=l} \bigotimes_{k=1}^r
\widetilde{H}_{l_k-1} (\, (\alpha_k,\, \hat{1}) \, ).
\end{eqnarray*}

(We can apply K\"unneth formula since we are working with chain
complexes of free modules over an hereditary ring $k$. The extra
$\textrm{Tor}$ terms in K\"unneth formula come from homology
groups of lower dimension which are null). If we define $m_k$ by
$|y_k|=m_k(t-1)+1$, the homology groups $\widetilde{H}_{l_k-1}(\,
(\alpha_k,\, \hat{1}) \, )$ are null for $l_k\neq m_k-1$.
Therefore, if $l$ is different from $\sum_{k=1}^r (m_k-1)$, we
have $\widetilde{H}_{l-1}(\, (x,\, y) \, )=0$. Since, the length
of maximal chains between $x$ and $y$ is equal to $d=\sum_{k=1}^r
(m_k-1)+1$ (\emph{cf.} proposition~\ref{grsubpos}), the homology
of the interval $[x,\, y]$ is concentrated in top dimension.

$(\Leftarrow)$ Conversely, if the poset $\Pi_\Po$ is
Cohen-Macaulay over the ring $k$, we have for any $m \geqslant 1$
and any $\alpha \in Min(\Pi_\Po(m(t-1)+1))$ that
$\widetilde{H}_{l-2}(\,(\alpha ,\, \hat{1}) \, )=0,$ if $l\ne m$.
Therefore, we get
$$H_{l}(\Pi_\Po(m(t-1)+1)) =
\bigoplus_{\alpha \in Min(\Pi_\Po(m(t-1)+1))} \widetilde{H}_{l-2}
(\, (\alpha,\, \hat{1}) \, )=0, $$ if $l\ne m$. And we conclude by
the theorem~\ref{KoszulPoset}.$\cqfd$
\end{deo}

\begin{rem}
To a graded poset, one can associate an algebra called
\emph{incidence algebra}. This incidence algebra is Koszul if and
only if the poset is Cohen-Macaulay (\emph{cf.} C. Cibils
\cite{Cibils}, P. Polo \cite{Polo} and D. Woodcock
\cite{Woodcock}). To a basic-set operad, one can define many
incidence algebras (one for each $n\in \mathbb{N}^*$ and each
$\alpha\in \Po_n([n])$). A direct corollary of this theorem and
the present work claims that a basic-set operad is Koszul if and
only if all its associated incidence algebras are Koszul.
\end{rem}

\section{Examples} \label{ExamplesPosets}

We introduce new partition type posets arising from operads. Using
the last theorem, we compute their homology in terms of the Koszul
dual cooperad.

\subsection{Homology of the partition poset
$\Pi$ and the operad $\C$} In this section, we mainly recall the
result of \cite{Fresse} which relates the homology of the
partition poset with the operad $\C$ of commutative algebras.\\

The operad $\C$ is a binary quadratic operad corresponding to
commutative algebras. Since there is only one multi-linear
operations with $n$ variables on a commutative algebra, the module
$\C_n$ has only one element.

\begin{dei}[Operad $\C$]
The operad $\C$ is defined by $\C_n:=\{e_n\}$, with $e_n$
invariant under the action of $\Sy_n$. Since $\C_n$ has only one
element the composition of the operad $\C$ is obvious. The operad
 $\C$ is a basic set operad.
\end{dei}

Since $\overline{e_n \times(x_1,\ldots,\, x_n)}= \overline{e_n
\times (x_{\sigma^{-1}(1)},\ldots,\, x_{\sigma^{-1}(n)})}$, we
have that $\C(I)$ has only one element which is corresponds to the
block $\{x_1,\ldots  ,\, x_n \}$. Therefore a $\C$-partition is a
classical partition of $[n]$ and the order between $\C$-partitions
is given by the refinement of partitions.

A. Bj\"orner showed in \cite{B} that the posets of partitions are
Cohen-Macaulay by defining an EL-labelling on them. Hence, the
homology of the partition poset is concentrated in top dimension
and given by
$$H_l(\Pi(n))\cong
\left\{
\begin{array}{ll}
\Li(n)^* \otimes
\textrm{sgn}_n& \textrm{if} \ l=n-1 \\
0 & \textrm{otherwise,}
\end{array}\right.$$
where $\Li(n)$ is the representation of $\Sn$ defined by the
operad $\Li$ of Lie algebras (\emph{cf.} \cite{Fresse, RW}).

\begin{rem}
P. Hanlon and M. Wachs have defined posets of partitions with
restricted block size in \cite{HW}. Each poset corresponds to an
operadic poset build on the operad $t-\C$ generated by a $t$-ary
commutative operation (\emph{cf.} \cite{Gnedbaye}). The same
EL-labelling as in the classical partition poset works here.
Therefore, we get that the operads $t-\C$ and their duals $t-\Li$
are Koszul over $\mathbb{Z}$ and over any field $k$.
\end{rem}

\subsection{Homology of the pointed partition poset
$\Pi_P$ and the operad $\Pe$} Following the same ideas, we
introduce the poset of pointed partitions coming from the operad
$\Pe$. We show that the homology of the pointed partition posets
is given by the Koszul dual cooperad of the operad $\Pe$ :
$\Pli^\vee$.

\subsubsection{\bf The operad $\bf \Pe$}
The operad $\Pe$, for \underline{perm}utation, has been introduced
by F. Chapoton in \cite{Chapoton}. \\

A \emph{$\Pe$-algebra} is an associative algebra such that
$x*y*z=x*z*y$ for any $x$, $y$ and $z$.

\begin{dei}[operad $\Pe$]
The operad $\Pe$ is the quadratic operad generated by a binary
operation $*$ and the following relations

\begin{enumerate}
\item \emph{(Associativity)} $(x*y)*z=x*(y*z)$ for any $x$, $y$
and $z$,

\item \emph{(Permutation)} $x*y*z=x*z*y$ for any $x$, $y$ and $z$.
\end{enumerate}
\end{dei}

We recall a result due to F. Chapoton. The algebraic operad $\Pe$
comes from the following set operad $\Pe_n:=\{e^n_1,\ldots,\,
e^n_n\}$ and, the composition $\mu$ is given by the formula :
\begin{eqnarray*}
\Pe_n\otimes\Pe_{i_1}\otimes
\cdots\otimes\Pe_{i_n} &\xrightarrow{\mu}& \Pe_{i_1+\cdots +i_n}\\
e^n_k\otimes e^{i_1}_{j_1} \otimes \cdots \otimes e^{i_n}_{j_n}
&\mapsto& e^{i_1+\cdots +i_n}_{i_1+\cdots +i_{k-1}+j_k}.
\end{eqnarray*}

(We denote by $\Pe$ the two operads.)
The image under the action of $\sigma\in\Sn$ of the vector $e^n_k$
is $e^n_{\sigma^{-1}(k)}$.

\begin{pro}
The operad $\Pe$ is a basic-set operad.
\end{pro}

\begin{deo}
The above formula of $\mu$ shows that every maps of the form
$\mu_{(\nu_1,\ldots,\, \nu_n)}$, where $\nu_l=e^{i_l}_{j_l}$, are
injective.$\cqfd$
\end{deo}

\subsubsection{\bf The posets $\bf \Pi_P$ of pointed partitions}
We describe the posets associated to the operad $\Pe$ in terms
of pointed partitions.\\

Since $\overline{e^n_i\times (x_1,\ldots,\, x_n)}=
\overline{e^n_{\sigma^{-1}(i)}\times (x_{\sigma^{-1}(1)},\ldots,\,
x_{\sigma^{-1}(n)})}$ in $\Pe(I)$, we choose to represent this
class of elements by $\{x_1,\ldots ,\, \bar{x_i},\ldots ,\, x_n
\}$, where the order between the $x_i$ does not care. The only
remaining information is which element is pointed. With this
identification, the set $\Pe_n([n])$ is equal to the set $\{\{1,\,
\ldots ,\, \bar{i},\, \ldots ,\, n \}\}_{1\leqslant i \leqslant
n}$, for any $n\in \mathbb{N}^*$.

\begin{dei}[Pointed partitions]
A \emph{pointed partition} $\{ B_1,\ldots,\, B_k \}$ of $[n]$ is a
partition on which one element of each block $B_i$ is emphasized.
\end{dei}

For instance, $\{ \{\bar{1},\, 3\},\, \{2,_, \bar{4}\}
\}$ is a pointed partition of $\{1,\, 2,\, 3,\, 4\}$.\\

The partial order on pointed partitions is a pointed variation of
the one for classical partitions. It is defined as follows.

\begin{dei}[Pointed partition poset, $\Pi_P$]
Let $\lambda$ and $\mu$ be two pointed partitions. We say that
$\mu$ is larger than $\lambda$ $(\lambda\leqslant \mu)$ if the
pointed integers of $\lambda$ belongs to the set of the pointed
integers of $\mu$ and if $\mu$ is a refinement of $\lambda$ as a
partition. This poset is called the \emph{pointed partition poset}
and denoted by $\Pi_P(n)$.
\end{dei}

For example, one has $\{\{\bar{1},\, 3\},\, \{2,\, \bar{4}\}\}
\leqslant \{\{\bar{1}\},\, \{\bar{3}\},\, \{2,\, \bar{4}\}\}$. The
largest element is $\{\{\bar{1}\},\ldots,$ $\{\bar{n}\}\}$. But
there are $n$ minimal elements of the type
$\{ \{1,\ldots,\,\bar{i},\ldots,\, n\} \}$. This poset is pure but not bounded.\\

We extend the action of the symmetric group in this case. Since
the image of a pointed integer under a permutation $\sigma\in\Sn$
gives a pointed integer, the symmetric group $\Sn$ acts on the
pointed partitions of $\{1,\ldots,\ n\}$.

\begin{pro}
The operadic poset $\Pi_{\Pe}$ associated to the operad $\Pe$ is
isomorphic to the poset of pointed partitions $\Pi_P$.
\end{pro}

\begin{deo}
In the case of the poset $\Pi_{\Pe}$, the order is defined by the
refinement of $\Pe$-partitions. With the identification
$\Pe_n([n])=\{\{1,\, \ldots ,\, \bar{i},\, \ldots ,\, n \}\}_{1\leqslant i
\leqslant n}$,
it corresponds to the refinement of pointed partitions.$\cqfd$
\end{deo}

\subsubsection{\bf Homology of the poset $\bf \Pi_P$}
We can use the properties of the operad $\Pe$ to compute the
homology of the pointed partition poset.

\begin{dei}[PreLie algebra]
A \emph{PreLie algebra} is a $k$-module $L$ equipped with a binary
operation $\circ$ such that
$$(x\circ y) \circ z- x\circ (y \circ z) =(x\circ z) \circ y - x\circ (z \circ y).$$
\end{dei}

We denote by $\Pli$ the related operad.

\begin{pro}{\rm [CV]}
\label{permkoszul} The pointed partition posets are totally
semi-modular. Recall that it implies that these posets are
Cohen-Macaulay over $k$.
\end{pro}

\begin{thm}[Homology of the pointed partition poset, $\Pi_P$]
The homology of the pointed partition poset is the $\Sn$-module
$$H_l(\Pi_P(n))\cong
\left\{
\begin{array}{ll}
\Pli(n)^* \otimes
\textrm{sgn}_{\Sy_n}& \textrm{if} \ l=n-1 \\
0 & \textrm{otherwise.}
\end{array}\right.$$
\end{thm}

\begin{deo}
The proof is a direct consequence of Theorems \ref{KoszulPoset}
and \ref{permkoszul}. $\cqfd$
\end{deo}

\begin{rem}
The operad $\Pli$ is isomorphic to the operad $\mathcal{RT}$ of
\emph{Rooted Trees} \cite{ChapLiv}. Therefore, the $\Sn$-module
$\Pli(n)$ is isomorphic to the free $k$-module on rooted trees
whose vertices are labelled by $\{1,\ldots\, n\}$, with the
natural action of $\Sn$ on them.
\end{rem}

We refer the reader to the article \cite{CV} of F. Chapoton and
the author for a study of the properties of the pointed partition
posets.

\subsection{Homology of the multi-pointed partition poset $\Pi_{MP}$
and the operad $\Comtri$}

To compute the homology of the multi-pointed partition poset, we
define a commutative version of the operad $\Trias$, which we call
the $\Comtri$ operad. We study its Koszul dual operad, called the
$\Postlie$ operad. (We prove in the appendix that these operads
are Koszul.)

\subsubsection{\bf The operad $\bf \Comtri$}

\begin{dei}[Commutative trialgebra]
A \emph{commutative trialgebra} is a $k$-module $A$ equipped with
two binary operations $*$ and $\bullet$ such that $(A,\, *)$ is a
$\Pe$-algebra
$$(x*y)*z=x*(y*z)=x*(z*y),$$
$(A,\,\bullet )$ is a commutative algebra
$$\left\{
\begin{array}{ccc}
x\bullet y&=&y \bullet x,\\
(x\bullet y)\bullet z&=&x\bullet (y\bullet z),
\end{array}
\right. $$ and the two operations $*$ and $\bullet$ must verify
the following compatibility relations
$$\left\{
\begin{array}{ccc}
x*(y\bullet z)&=&x*(y*z),\\
(x\bullet y)*z&=&x\bullet (y*z).
\end{array}
\right. $$
\end{dei}

\begin{dei}[$\Comtri$ operad]
We denote by $\Comtri$, the operad coding the commutative
trialgebras.
\end{dei}

The operad $\Comtri$ is an operad generated by two operations $*$
and $\bullet$ and the relations defined before. As a consequence
of Theorem~\ref{FreeComTri}, giving the free commutative
trialgebra, we have the following complete description of this
operad.

\begin{thm}
\label{ComTriSet}
 For any $n\leqslant1$, the $\Sn$-module
$\Comtri(n)$ is isomorphic to the free $k$-module on $\Comtri_n=\{
e^n_J \ / \ J\subset [n],\ J\ne \emptyset\}$, where the action of
$\, \Sn$ is induced by the natural action on the subsets $J$ of
$[n]$.

The composition of the operad $\Comtri$ is given by the formula
\begin{eqnarray*}
\Comtri_k\otimes \Comtri_{i_1}\otimes \cdots \otimes \Comtri_{i_k}
&\to& \Comtri_{i_1+\cdots+i_k=n} \\
e^k_J\otimes e^{i_1}_{J_1}\otimes \cdots \otimes e^{i_k}_{J_k}
&\mapsto& e^n_{\bar{J}},
\end{eqnarray*}
where $\bar{J}=\bigcup_{j\in J} (i_1+\cdots +i_{j-1})+J_j$. (The
subset of $[n]$ denoted by $(i_1+\cdots +i_{j-1})+J_{j}$ is
defined by $\bigcup_{a\in J_j} \{i_1+\cdots +i_{j-1}+a\}$.)
\end{thm}

\begin{deo}
Let $\Po$ be the $\Sy$-module generated by the sets $\Comtri_n$
($\Po(n) = k[\Comtri_n]$). We consider the Schur functor
associated to $\Po$. For any $k$-module $V$, we have
$\mathcal{S}_\Po(V)=\bigoplus_{n\geqslant 1}
k[\beta_n]\otimes_{\Sy_n} V^{\otimes n}$. We denote the element
 $e_J^n\otimes (x_1, \ldots ,\, x_n)$ of
$k[\Comtri_n]\otimes V^{\otimes n}$ by $x_{j_1}\ldots x_{j_k}
\otimes x_{l_1} \ldots x_{l_m}$, where $J=\{j_1,\ldots ,\, j_k \}$
and $[n]-J=\{l_1,\ldots ,\, l_m \}$. Therefore, we get
\begin{eqnarray*}
\mathcal{S}_\Po(V)=\bigoplus_{n\geqslant 1} \bigoplus_{j=1}^n
\bigoplus_{|J|=j} k[e^n_J]\otimes_{\Sy_n} V^{\otimes n} \cong
\bigoplus_{n\geqslant 1} \bigoplus_{j=1}^n \bar{S}_j(V)\otimes
S_{n-j}(V)= \bar{S}(V)\otimes S(V).
\end{eqnarray*}
Since $\mathcal{S}_\Po(V)$ is equal to $\mathcal{S}_{\Comtri}(V)$
for any $k$-module $V$, the $\Sy$-module $\Po$ is isomorphic  to
$\Comtri$. With the following identifications
$$ e^2_{\{1\}} \longleftrightarrow * \quad , \quad
e^2_{\{2\}} \longleftrightarrow *^{(1,\, 2)} \quad  , \quad
e^2_{\{1,\, 2\}} \longleftrightarrow \bullet ,$$ one can easily
see that $\mathcal{S}_\Po(V)$ is the free commutative trialgebra
defined in the previous theorem. The following diagram is
commutative.
$$\xymatrix{\mathcal{S}_\Po\circ \mathcal{S}_\Po \ar[r] \ar[d]^{\mu_\Po}
& \mathcal{S}_{\Comtri}\circ \mathcal{S}_{\Comtri}  \ar[d]^{\mu_{\Comtri}}\\
\mathcal{S}_\Po \ar[r] &   \mathcal{S}_{\Comtri}  }$$ We conclude
that the operads $\Po$ and $\Comtri$ are isomorphic. $\cqfd$
\end{deo}

\begin{cor}
The operad $\Comtri$ is a basic-set operad.
\end{cor}

\begin{deo}
It remains to prove that on the $\Comtri_n$ the maps
$\mu_{(\nu_1,\ldots,\, \nu_k)}$ are injective, which is straightforward.$\cqfd$
\end{deo}

\begin{rem}
The elements of $\Comtri(n)$ can be indexed by the cells of the
simplex of dimension $n-1$ : $\Delta^{n-1}$. The action of $\Sn$
on $\Comtri(n)$ corresponds to the action of $S_n$ on the faces of
$\Delta^{n-1}$.
\end{rem}

\subsubsection{\bf The posets $\bf \Pi_{MP}$ of multi-pointed partitions}
We describe the posets $\Pi_{\Comtri}$ in terms of multi-pointed partitions.\\

Following the same idea as in the previous case, we choose to
represent the element $\overline{e^n_J\times (x_1,\ldots ,\,
x_n)}$ of $\Comtri_n([n])$ by $\{x_1,\ldots ,\, \overline{x_{j_1}}
,\ldots ,\, \overline{x_{j_k}} ,\ldots ,\, x_n \}$, where
$J=\{j_1, \ldots ,\,j_k \}$.

\begin{dei}[Multi-pointed partition]
A \emph{multi-pointed partition} $\{ B_1,\ldots,\, B_k \}$ of
$[n]$ is a partition on which at least one element of each block
$B_i$ is emphasized.
\end{dei}

\begin{dei}[Multi-pointed partition poset, $\Pi_{MP}$]
Let $\lambda$ and $\mu$ be two multi-pointed partitions. We say
that $\mu$ is larger than $\lambda$ $(\lambda\leqslant \mu)$ if
 $\mu$ is a refinement of
$\lambda$ as a partition and if for each block $B$ of $\mu$, all
the pointed integers of $B$ remain pointed in $\lambda$ or all
become unpointed. This poset is called the \emph{multi-pointed
partition poset} and denoted by $\Pi_{MP}(n)$.
\end{dei}

The action of the symmetric group is defined in a similar way and
preserves the partial order $\leqslant$. With this identification,
one can see that the refinement of $\Comtri$-partitions
corresponds to the refinement of multi-pointed partitions.
Therefore, we have the following proposition.

\begin{pro}
The operadic poset $\Pi_{\Comtri}$  is isomorphic to the poset of
multi-pointed pointed partitions $\Pi_{MP}$.
\end{pro}

\subsubsection{\bf Homology of the poset
$\bf \Pi_{MP}$}

\begin{dei}[PostLie algebra]
A \emph{PostLie algebra} is a $k$-module $L$ equipped with two
binary operations $\circ$ and $[\, ,\,]$ such that $(L,\,
[\,,\,])$ is a Lie algebra
$$\left\{
\begin{array}{l}
\lbrack x,\, y \rbrack=-\lbrack y ,\, x\rbrack \\
\lbrack x,\, \lbrack y,\, z\rbrack \rbrack + \lbrack y,\, \lbrack
z,\, x\rbrack \rbrack + \lbrack z,\, \lbrack x,\, y \rbrack
\rbrack =0\end{array} \right.$$ and such that the two operations
$\circ$ and $[\, ,\, ]$ verify the following compatibility
relations
$$\left\{
\begin{array}{l}
(x\circ y)\circ z - x\circ (y\circ z) - (x\circ z)\circ y + x\circ
(z\circ y) = x\circ \lbrack y,\, z\rbrack, \\
\lbrack x,\, y \rbrack \circ z = \lbrack x\circ z ,\, y \rbrack
+\lbrack x,\,  y\circ z\rbrack
\end{array}
\right. $$
\end{dei}

\begin{dei}[$\Postlie$ operad]
We denote by $\Postlie$, the operad coding the PostLie algebras.
\end{dei}

Theorem~\ref{ComtriKoszulPostlie} proves that the operad
$\Postlie$ is the Koszul dual of $\Comtri$. In \cite{CV}, we have
shown that the posets of multi-pointed partitions are totally
semi-modular. Therefore, they are Cohen-Macaulay and we have the
following theorem.

\begin{thm}[Homology of the multi-pointed partition poset, $\Pi_{MP}$]
The homology of the multi-pointed partition poset $\Pi_{MP}$ with
coefficient in $k$ is the $\Sn$-module
$$H_l(\Pi_{MP}(n))\cong
\left\{
\begin{array}{ll}
(\Li\circ \Mag)^*(n)\otimes sgn_{\Sy_n}
& \textrm{if} \ l=n-1 \\
0 & \textrm{otherwise.}
\end{array}\right.$$
\end{thm}

\begin{rem}
We refer the reader to \cite{CV} for further properties on the
multi-pointed partition posets.
\end{rem}

\subsection{Homology of the ordered partition poset $\Pi_O$
and the operad $\A$} The posets associated to the operad $\A$ are
the posets of ordered partitions (or bracketing). The homology of
these posets is given by the Koszul dual cooperad of the operad
$\A$ : $\A^\vee$.

\subsubsection{\bf The operad $\A$}

The operad $\A$ corresponds to the operad of associative algebras.
This operad is generated by one non-symmetric binary operation and
the associativity relation. One has a complete description of it.
Up to permutation, there exists only one associative operation on
$n$ variables. Therefore, we have $\A(n)=k\lbrack \Sn \rbrack$,
the regular representation of $\Sn$.\\

The operad $\A$ comes from the set operad $\A_n=\Sy_n$ with the
following composition maps
\begin{eqnarray*}
\Sy_n\times\Sy_{ i_1}\times
\cdots\times\Sy_{i_n} &\xrightarrow{\mu}& \Sy_{i_1+\cdots +i_n}\\
\pi\times (\tau_1, \ldots,\, \tau_n) &\mapsto&
[\tau_{\pi^{-1}(1)},\ldots,\, \tau_{\pi^{-1}(n)}],
\end{eqnarray*}
where $[\tau_{\pi^{-1}(1)},\ldots,\, \tau_{\pi^{-1}(n)}]$ is the
block permutation of $\Sy_{i_1+\cdots +i_n}$ induced by

$\tau_{\pi^{-1}(1)},\ldots,\, \tau_{\pi^{-1}(n)}$.\\

The image under the action of $\sigma \in \Sy_n$ of the vector
$\tau\in \beta_n=\Sy_n$ is $\tau \sigma$.

\begin{pro}
The operad $\A$ is a basic-set operad.
\end{pro}

\begin{deo}
 The composition maps $\mu_{(\tau_1,
\ldots,\, \tau_n)}$ are clearly injective.$\cqfd$
\end{deo}

\subsubsection{\bf The posets $\Pi_O$ of ordered partitions}

For any $n\in \mathbb{N}^*$, we choose to represent the element
$$\overline{\sigma\times (x_1,\ldots ,\, x_n)}=
\overline{id \times (x_{\sigma(1)},\ldots ,\, x_{\sigma(n)})}$$ of $\A_n([n])$ by
$ (x_{\sigma(1)},\ldots ,\, x_{\sigma(n)})$.

\begin{dei}[Ordered partitions]
Any set $\{ B_1,\ldots,\, B_r\}$ of ordered sets of elements of
$[n]$ that gives a partition when forgetting the order between the
integers of each block is called an \emph{ordered partition}.
\end{dei}

For instance $\{(3,\, 1),\,(2,\, 4)\}$ and
$\{(1),\,(3),\,(2,\,4)\}$
are two ordered partitions of $\,$  $\{1,2,\, 3,\,4\}$.\\

A partial order on the set of ordered partitions is defined as
follows.

\begin{dei}[Ordered partition poset, $\Pi_O$]
Let $\lambda=\{ B_1,\ldots,\, B_r\}$ and $\mu=\{ C_1,\ldots$, $
C_s\}$ be two ordered partitions of $[n]$. We say that $\mu$ is
larger than $\lambda$ $(\lambda\leqslant \mu)$ if, for every
$B_i$, there exists an ordered sequence of $C_j$ such that
$B_i=(C_{j_1}, \ldots, \, C_{j_i})$.
 We call this poset the
\emph{ordered partition poset} and we denote it by $\Pi_O$.
\end{dei}

In the previous example, $\{(3,\, 1),\,(2,\, 4)\} \leqslant
\{(1),\,(3),\,(2,\,4)\}$. The minimal elements are of the form $\{
(\tau(1),\ldots,\, \tau(n))\}$, for $\tau \in \Sn$, and the
maximal
one is $\{ (1),\ldots ,\,(n)\}$. This poset is pure.\\

The ordered partition poset is equipped with the natural action of
the symmetric group $\Sn$ which is compatible with the partial
order. Observe that the action of $\Sn$ on the set
$\Delta(\Pi_O(n))$ is free. All the intervals $[\alpha,\,
\hat{1}]$, with $\alpha\in Min(\Pi_O(n))$, are isomorphic and the
isomorphism between two of them is given by the action of an
element of $\Sn$.

\begin{pro}
The operadic poset $\Pi_{\A}$ associated to the operad $\A$ is
isomorphic to the poset of ordered partitions $\Pi_O$.
\end{pro}

\subsubsection{\bf Homology of the poset $\Pi_O$}

\begin{pro}
\label{Askoszul} Every interval of the from $[\alpha,\, \hat{1}]$
with $\alpha\in Min(\Pi_O(n))$ is isomorphic to the boolean
lattice of subsets of $[n-1]$. Therefore, the intervals
$[\alpha,\, \hat{1}]$ with $\alpha\in Min(\Pi_O(n))$ are
Cohen-Macaulay over $k$.
\end{pro}

\begin{deo}
Since these intervals are all isomorphic, it is enough to prove it
for one of them. Let $\alpha$ be the ordered partition $(1,\, 2,\,
\ldots ,\, n)$ coming from the identity in $\Sn=\A_n$. To an
ordered partition $(1,\ldots,\, i_1),\, (i_1+1,\ldots ,\,
i_2),\ldots ,\,(i_k+1,\ldots,\, n)$ of $[\alpha,\, \hat{1}]$ with
$i_1<i_2<\ldots<i_k$, one can associate the subset
$\{i_1,\ldots,\,i_k \}$ of $[n-1]$. This map defines an
isomorphism of posets.
\end{deo}

\begin{thm}[Homology of the ordered partition poset, $\Pi_O$]
The homology of the ordered partition poset $\Pi_O$ is the
$\Sn$-module
$$H_l(\Pi_O(n))\cong
\left\{
\begin{array}{ll}
k[\Sn]& \textrm{if} \ l=n-1 \\
0 & \textrm{otherwise.}
\end{array}\right.$$
\end{thm}

\begin{rem}
Proposition \ref{Askoszul} proves that the operad $\A$ is Koszul
over $k$.
\end{rem}

\subsection{Homology of the pointed ordered partition poset
 $\Pi_{PO}$ and the operad $\D$}

From the operad $\D$ of associative dialgebras, we construct the
pointed ordered partition posets. The homology of the pointed
ordered partition poset $\Pi_{PO}$ is given by the Koszul dual
cooperad of the operad $\D$ : $\Dend^\vee$.

\subsubsection{\bf The operad $\D$}

J.-L. Loday has defined the notion of \emph{associative dialgebra}
(\emph{cf.} \cite{Loday}) in the following way. A
\emph{associative dialgebra} is a vector space $D$ equipped with
two binary operations $\dashv$ and $\vdash$ such that
$$
\left\{
\begin{array}{l}
 (x \dashv y ) \dashv z = x\dashv (y \vdash z)\\
 (x\dashv y )\dashv z = x \dashv(y \dashv z)\\
 (x\vdash y ) \dashv z = x \vdash(y \dashv z)\\
 (x \dashv y )\vdash z = x \vdash(y\vdash z)\\
 (x \vdash y )\vdash z = x\vdash (y \vdash z)
\end{array}
\right.
$$

The associated binary quadratic operad is denoted by $\D$.\\

We have seen that the operad $\Pe$ is a pointed version of $\C$.
The operad $\D$ is a pointed version of $\A$ mainly because
$\D=\Pe\otimes\A$ (\emph{cf.} \cite{Chapoton}). Each operation on
$n$ variables in $\D(n)$, is an associative operation where one
input is emphasized. Therefore, the $\Sn$-module $\D(n)$ is
isomorphic to $n$ copies of $\A(n)$:
$$\D(n)=k^n \otimes k[\Sn].$$

This operad comes from the following set operad
$\D_n=\{e^n_k\}_{
1\leqslant  k \leqslant n}\times \Sy_n$.

\begin{pro}
The operad $\D$ is a basic-set operad.
\end{pro}

\begin{deo}
The proof is a direct consequence of the case $\Pe$ and $\A$.
$\cqfd$
\end{deo}

\subsubsection{\bf The posets $\Pi_{PO}$ of pointed ordered
partitions}

The set $\D_n([n])$ corresponds to the set $\{( \tau(1), \ldots
,\, \overline{\tau(i)},\ldots  ,\, \tau(n) )\}_{\tau \in \Sy_n,\,
1\leqslant i\leqslant n}$.

\begin{dei}[Pointed ordered partition]
A \emph{pointed ordered partition} is an ordered partition
$\{B_1,\ldots,\, B_k\}$ of $[n]$ where, for each $1 \leqslant i
\leqslant k$, an element of the sequence $B_i$ is emphasized.
\end{dei}

The partial order on pointed ordered partitions is a pointed
variation of the one on ordered partitions.

\begin{dei}[Pointed ordered partition poset, $\Pi_{PO}$]
Let $\lambda$ and $\mu$ be two pointed ordered partitions. We say
that $\mu$ is larger than $\lambda$ if the pointed integers of
$\lambda$ belongs to the set of the pointed integers of $\mu$ and
if $\mu$ is a larger then $\lambda$ as ordered partitions. This
poset is called the \emph{pointed ordered partition poset} and
denoted by $\Pi_{PO}$.
\end{dei}

For example, we have $\{(3,\, \bar{1}),\, (2,\, \bar{4})\}
\leqslant \{(\bar{1}) ,\, (\bar{3}),\, (2,\, \bar{4})\}$. Any
element of the form $\{(\tau(1),\ldots,\, \overline{\tau(i)},
\ldots,\, \tau(n))\}$ is a minimal element in this poset, and
$\{(\bar{1}),\ldots
,\, (\bar{n})\}$ is the maximal element. Once again, this poset is pure.\\

The action of the symmetric group $\Sn$ is defined in a similar
way and is compatible with the partial order. Once again, the
action of $\Sn$ on the set $\Delta(\Pi_{PO}(n))$ is free. All the
intervals $[\alpha,\, \hat{1}]$, with $\alpha\in
Min(\Pi_{PO}(n))$, are isomorphic and the isomorphism between two
of them is given by the action of an element of $\Sn$.

\begin{pro}
The operadic poset $\Pi_{\D}$ associated to the operad $\D$ is
isomorphic to the poset of pointed ordered partitions.
\end{pro}

\subsubsection{\bf Homology of the poset $\Pi_{PO}$}

The Koszul dual of the operad $\D$ is the operad $\Dend$ of
dendriform algebras (\emph{cf.} \cite{Loday}). Up to permutations,
the linear operations on $n$ variables of a dendriform algebra can
be indexed by the planar binary trees with $n$ vertices $Y_n$. One
has $\Dend(n)\simeq k[Y_n]\otimes k\lbrack \Sn \rbrack$.

\begin{pro}
Every interval $[\alpha,\, \hat{1}]$ of  $\Pi_{PO}(n)$ with
$\alpha\in Min(\Pi_{PO}(n))$ is totally semi-modular and therefore
Cohen-Macaulay over $k$.
\end{pro}

\begin{deo}
The proof is the same as the one for pointed partitions in
\cite{CV}.
\end{deo}

\begin{thm}[Homology of the pointed ordered poset, $\Pi_{PO}$]
The homology of the pointed ordered partition poset $\Pi_{PO}$ is
the $\Sn$-module
$$H_l(\Pi_{PO}(n))\cong
\left\{
\begin{array}{ll}
k[Y_n]\otimes k[\Sn]
& \textrm{if} \ l=n-1 \\
0 & \textrm{otherwise.}
\end{array}\right.$$
\end{thm}

\begin{rem}
Proposition \ref{Askoszul} proves that the operads $\D$ and
$\Dend$ are Koszul over $k$.
\end{rem}

\subsection{Homology of the multi-pointed ordered partition poset
 $\Pi_{MPO}$ and the operad $\Trias$}

The multi-pointed ordered partition posets $\Pi_{MPO}$ arise from
the operad $\Trias$ of associative trialgebras. Therefore the
homology of the multi-pointed ordered partition poset is given by
the Koszul dual cooperad of $\Trias$, namely $\Tridend^\vee$. \\

The notion of an \emph{associative trialgebra} has been defined by
J.-L. Loday and M.O. Ronco in \cite{LodayRonco}. Following the
same methods as before, one can prove that the operad $\Trias$ is
a basic-set operad. The corresponding partition poset is
isomorphic to the poset of multi-pointed ordered partitions.

\begin{dei}[Multipointed ordered partition]
A \emph{multipointed ordered partition} is an ordered partition
$\{ B_1,\ldots,\, B_k \}$ of $[n]$ where, for each $i$, at least
one element of the sequence $B_i$ is emphasized.
\end{dei}

\begin{dei}[Multipointed ordered partition poset, $\Pi_{MPO}$]
Let $\lambda$ and $\mu$ be two multipointed ordered partitions. We
say that $\mu$ is larger than $\lambda$ if the pointed integers of
$\lambda$ belongs to the set of the pointed integers of $\mu$ and
if $\mu$ is a larger than $\lambda$ as ordered partitions. This
poset is called the \emph{multipointed ordered partition poset}
and denoted by $\Pi_{MPO}$.
\end{dei}

\begin{pro}
Every interval $[\alpha,\, \hat{1}]$ of  $\Pi_{MPO}(n)$ with
$\alpha\in Min(\Pi_{MPO}(n))$ is totally semi-modular and
therefore Cohen-Macaulay over $k$.
\end{pro}

The $\Sy_n$-modules of the operad $\Tridend$ are given by
$\Tridend(n)=k[T_n]\otimes k[\Sy_n]$, where $T_n$ denotes the set
of planar trees with $n$ leaves. (It corresponds to the cells of
the associahedra.)

\begin{thm}[Homology of the multi-pointed ordered partition poset, $\Pi_{MPO}$]
The homology of the multi-pointed ordered partition poset
$\Pi_{MPO}$ is the $\Sn$-module
$$H_l(\Pi_{MPO}(n))\cong
\left\{
\begin{array}{ll}
k[T_n]\otimes k[\Sn]
& \textrm{if} \ l=n-1 \\
0 & \textrm{otherwise.}
\end{array}\right.$$
\end{thm}

\subsection{Comparison between the various partition posets}

We sum up the various cases in the following table. \\

\begin{tabular}{|c|c|c|c|}
\hline
\textbf{Partition poset} $\bf \Pi_{\Po}$ & \textbf{Operad} $\bf \Po$ & $\bf \Po^!$ &
$\bf H_{n-1}(\Pi_\Po(n))=\Po^{\ac}(n)$ \\
\hline
classical $\Pi$ & $\C$ & $\Li$ & $\Li^*(n) \otimes sgn_{\Sy_n}$  \\
\hline $t$-restricted block size &  $t-\C$ & $t-\Li$ &  \\
\hline pointed  $\Pi_P$ & $\Pe$ & $\Pli$ & $\mathcal{RT}^*(n)\otimes sgn_{\Sy_n}$ \\
\hline multi-pointed $\Pi_{MP}$ & $\Comtri$ & $\Postlie$ & $(\Li
\circ
\Mag)^*(n)\otimes sgn_{\Sy_n} $ \\
\hline ordered  $\Pi_O$ & $\A$  & $\A$  &  $k[\Sy_n]$ \\
\hline pointed ordered $\Pi_{PO}$ & $\D$ & $\Dend$ & $k[Y_n]\otimes k[\Sy_n]$\\
\hline multi-pointed ordered $\Pi_{MPO}$ & $\Trias$ & $\Tridend$ & $k[T_n] \otimes k[\Sy_n]$\\
\hline
\end{tabular}

$ $\\

We can define inclusion and forgetful maps linking the six cases
of partition posets defined above.

\begin{dei}[Forgetful maps $F_O$ and $F_P$]
The forgetful map \emph{$F_O$} is obtained from an ordered type
partition by forgetting the internal order inside each block.

The forgetful map \emph{$F_P$} unmark all emphasized integers of a
pointed partition (resp. a pointed ordered partition) in order to
get a partition (resp. an ordered partition).
\end{dei}

Since a pointed type partition is a multipointed type partition,
this defines an injective map from pointed type partitions to
multipointed type partitions.\\

These maps induce morphisms of presimplicial sets that commute
with the action of the symmetric group $\Sn$. We sum up the six
cases studied above in the following diagram.

$$\xymatrix{\Pi_{MPO} \ar[r]^{F_O} & \Pi_{MP} \\
\Pi_{PO} \ar[u] \ar[d]^{F_P}  \ar[r]_{F_O}&  \Pi_P \ar[u] \ar[d]^{F_P}\\
\Pi_O \ar[r]_{F_O} & \Pi  }$$

 The comparison diagram between the
various posets corresponds to the following diagram of operads.
The next diagram gives the homology of the various posets.

$$\xymatrix{ \Trias \ar[d]  \ar[r] &   \ar[d]\Comtri & & \Tridend  & \ar[l] \Postlie \\
\D  \ar[r] & \Pe & & \Dend  \ar[u] \ar[d]  & \Pli \ar[l] \ar[u] \ar[d] \\
\A \ar[u]  \ar[r] &  \ar[u]\C & & \A &  \ar[l]\Li}$$

\newpage

\appendix

\section*{ APPENDIX : Koszul duality of the operads $\Comtri$ and $\Postlie$}
We define a commutative version of the operad $\Trias$, which we
call the $\Comtri$ operad. We study its Koszul dual operad, which
we call the $\Postlie$ operad. We prove that the operads $\Comtri$
and $\Postlie$ are Koszul operads.\\

Let $k$ be a field of characteristic $0$.

\begin{flushleft}
1. $ $ {\bf   $\bf \Comtri$ operad}
\end{flushleft}
$ $

\begin{dei}[Commutative trialgebra]
A \emph{commutative trialgebra} is a $k$-module $A$ equipped with
two binary operations $*$ and $\bullet$ such that $(A,\, *)$ is a
$\Pe$-algebra
$$(x*y)*z=x*(y*z)=x*(z*y),$$
$(A,\,\bullet )$ is a commutative algebra
$$\left\{
\begin{array}{ccc}
x\bullet y&=&y \bullet x,\\
(x\bullet y)\bullet z&=&x\bullet (y\bullet z),
\end{array}
\right. $$ and the two operations $*$ and $\bullet$ must verify
the following compatibility relations
$$\left\{
\begin{array}{ccc}
x*(y\bullet z)&=&x*(y*z),\\
(x\bullet y)*z&=&x\bullet (y*z).
\end{array}
\right. $$
\end{dei}

\begin{dei}[$\Comtri$ operad]
We denote by $\Comtri$, the operad coding the commutative
trialgebras.
\end{dei}

\begin{rem}
The operad $\Comtri$ does not fall into the construction of M.
Markl (\emph{cf.} \cite{Markl}) called \emph{distributive laws}
because of the first compatibility relation. Therefore, we will
have to use other methods to show that this operad is Koszul.
\end{rem}

\begin{thm}[Free commutative trialgebra]
\label{FreeComTri} The free commutative trialgebra algebra on a
module $V$, denoted by $\Comtri(V)$ is given by the module
$\bar{S}(V)\otimes S(V)$ equipped with the following operations
\begin{eqnarray*}
(x_1\ldots x_k\otimes y_1\ldots y_l)*(x'_1\ldots x'_m\otimes
y'_1\ldots y'_n)=x_1\ldots x_k \otimes y_1\ldots y_l \, x'_1\ldots
x'_m\,  y'_1\ldots y'_n, \\
(x_1\ldots x_k\otimes y_1\ldots y_l)\bullet(x'_1\ldots x'_m\otimes
y'_1\ldots y'_n)=x_1\ldots x_k\,  x'_1\ldots x'_m\otimes y_1\ldots
y_l\,  y'_1\ldots y'_n.
\end{eqnarray*}
\end{thm}

\begin{deo}
It is easy to check that the $k$-module $\Comtri(V)$ equipped with
the two operations $*$ and $\bullet$ is a commutative trialgebra.

Let $(T,\, *,\,\bullet) $ be a commutative trialgebra. The
inclusion $V \hookrightarrow V\otimes 1$ is denoted by $i$.
Associated to each morphism of $k$-modules $f \, : \, V \to T$, we
will show that there exists is a unique morphism of commutative
trialgebras $\tilde{f}\, : \, \bar{S}(V)\otimes S(V)\to T $ such
that the following diagram is commutative

$$\xymatrix{V \ar[r]^(0.3){i} \ar[dr]_{f}& \bar{S}(V)\otimes S(V) \ar[d]^{\tilde{f}} \\
 & T.} $$

We define $\tilde{f}$ by the formula
\begin{eqnarray*}
\tilde{f}(X\otimes Y)=\tilde{f}(x_1\ldots x_k\otimes y_1\ldots
y_l) &:=& \big( \tilde{f}(x_1)\bullet \cdots \bullet
\tilde{f}(x_k)\big)* \big(
\tilde{f}(y_1)* \cdots * \tilde{f}(y_l)\big)\\
\tilde{f}(x_1\ldots x_k\otimes 1) &:=& \tilde{f}(x_1)\bullet
\cdots \bullet \tilde{f}(x_k).
\end{eqnarray*}

The map $\tilde{f}$ is a morphism of commutative trialgebras. We
have

\begin{eqnarray*}
\tilde{f}(X\otimes Y * X'\otimes Y') &=& \tilde{f}(x_1\ldots x_k
\otimes y_1\ldots y_l\,  x'_1\ldots x'_m\, y'_1\ldots y'_n) \\
&=& \big( \tilde{f}(x_1)\bullet \cdots \bullet \tilde{f}(x_k)
\big)* \big( \tilde{f}(y_1)* \cdots * \tilde{f}(y_l)
* \\
& & \tilde{f}(x'_1)* \cdots * \tilde{f}(x'_n)
* \tilde{f}(y'_1)* \cdots * \tilde{f}(y'_m) \big) \\
&=& \big( \tilde{f}(x_1)\bullet \cdots \bullet \tilde{f}(x_k)\big)
* \big( (\tilde{f}(y_1) * \cdots * \tilde{f}(y_l))
* \\
& & (\tilde{f}(x'_1) \bullet \cdots \bullet \tilde{f}(x'_n))
* (\tilde{f}(y'_1)* \cdots * \tilde{f}(y'_m)) \big) \\
&=& \big( (\tilde{f}(x_1)\bullet \cdots \bullet \tilde{f}(x_k))
*  (\tilde{f}(y_1) * \cdots * \tilde{f}(y_l)) \big)
* \\
& & \big( (\tilde{f}(x'_1) \bullet \cdots \bullet \tilde{f}(x'_n))
* (\tilde{f}(y'_1)* \cdots * \tilde{f}(y'_m)) \big) \\
&=& \tilde{f}(X\otimes Y) * \tilde{f}(X'\otimes Y'),
\end{eqnarray*}

and

\begin{eqnarray*}
\tilde{f}(X\otimes Y \bullet X'\otimes Y') &=& \tilde{f}(x_1\ldots
x_k \,  x'_1\ldots x'_m
\otimes y_1\ldots y_l\, y'_1\ldots y'_n) \\
&=& \big( \tilde{f}(x_1)\bullet \cdots \bullet \tilde{f}(x_k)
\bullet \tilde{f}(x'_1)\bullet \cdots \bullet
\tilde{f}(x'_m)\big)* \\
& & \big( \tilde{f}(y_1)* \cdots * \tilde{f}(y_l)
* \tilde{f}(y'_1)* \cdots *
\tilde{f}(y'_n) \big).
\end{eqnarray*}

In any commutative trialgebra, one has

\begin{eqnarray*}
(a\bullet a')*(b * b')&=& (a'\bullet a)*(b * b') = (a'\bullet
(a*b))
* b'= ((a*b)\bullet a')*b'\\
&=& (a*b)\bullet (a'*b').
\end{eqnarray*}

Therefore, we get

\begin{eqnarray*}
\tilde{f}(X\otimes Y \bullet X'\otimes Y') &=& \big(
(\tilde{f}(x_1)\bullet \cdots \bullet \tilde{f}(x_k)) *
(\tilde{f}(y_1)* \cdots * \tilde{f}(y_l))\big) \bullet \\
& & \big( (\tilde{f}(x'_1)\bullet \cdots \bullet \tilde{f}(x'_m))
*
(\tilde{f}(y'_1)* \cdots * \tilde{f}(y'_n))\big)\\
&=&\tilde{f}(X\otimes Y) \bullet \tilde{f}(X'\otimes Y').
\end{eqnarray*}

Let $g \, : \, \Comtri(V) \to T$ be a morphism of commutative
trialgebras such that $g \circ i = f$. Since a tensor $x_1\ldots
x_k\otimes y_1\ldots y_l$ in $\bar{S}(V)\otimes S(V)$ is equal  to
$(x_1\otimes 1 \bullet \cdots \bullet x_k\otimes 1 )*(y_1\otimes 1
* \cdots * y_l\otimes 1)$, we have
\begin{eqnarray*}
g(X\otimes Y) &=& g\big( (x_1\otimes 1 \bullet \cdots \bullet
x_k\otimes 1)*(y_1\otimes 1 * \cdots * y_l\otimes 1) \big) \\
&=& \big( g(x_1\otimes 1) \bullet \cdots \bullet
g(x_k\otimes 1)\big) * \big( g(y_1\otimes 1) * \cdots * g(y_l\otimes 1) \big) \\
&=&\big( f(x_1) \bullet \cdots \bullet f(x_k)\big)
* \big( f(y_1) * \cdots * f(y_l) \big) \\
&=& \tilde{f}(X\otimes Y).
\end{eqnarray*}
$\cqfd$
\end{deo}

We recall from theorem~\ref{ComTriSet} that the $\Sy_n$-modules
$\Comtri(n)$ are free $k$-modules.

\begin{rem}
There exists a differential graded variation of the operad $\Comtri$ defined
by F. Chapoton in \cite{ChapotonLambda}.
\end{rem}

\newpage
\begin{flushleft}
2. $ $ {\bf  $\bf \Postlie$ operad }
\end{flushleft}

\begin{dei}[PostLie algebra]
A \emph{PostLie algebra} is a $k$-module $L$ equipped with two
binary operations $\circ$ and $[\, ,\,]$ such that $(L,\,
[\,,\,])$ is a Lie algebra
$$\left\{
\begin{array}{l}
\lbrack x,\, y \rbrack=-\lbrack y ,\, x\rbrack \\
\lbrack x,\, \lbrack y,\, z\rbrack \rbrack + \lbrack y,\, \lbrack
z,\, x\rbrack \rbrack + \lbrack z,\, \lbrack x,\, y \rbrack
\rbrack =0\end{array} \right.$$ and such that the two operations
$\circ$ and $[\, ,\, ]$ verify the following compatibility
relations
$$\left\{
\begin{array}{l}
(x\circ y)\circ z - x\circ (y\circ z) - (x\circ z)\circ y + x\circ
(z\circ y) = x\circ \lbrack y,\, z\rbrack, \\
\lbrack x,\, y \rbrack \circ z = \lbrack x\circ z ,\, y \rbrack
+\lbrack x,\,  y\circ z\rbrack
\end{array}
\right. $$
\end{dei}

\begin{dei}[$\Postlie$ operad]
We denote by $\Postlie$, the operad coding the PostLie algebras.
\end{dei}

\begin{rem}
A Lie algebra $(L,\, [,])$ equipped an extra operation $\circ=0$
is a PostLie algebra. A PostLie algebra $(L,\, \circ,\, [\, , \,
])$ such that  $[\, ,\,]=0$ is a PreLie algebra for the product
$\circ$. Therefore the operad $\Postlie$ is an ``extension'' of
the operad $\Pli$ by the operad $\Li$.
\end{rem}

We will compute the free PostLie algebra. It is expressed with the
free magmatic algebra $(\Mag(V),\, \diamond)$. A \emph{magmatic
algebra} is a $k$-module $A$ equipped with a binary operation.
Therefore, the free magmatic algebra on a $k$-module $V$ is
defined on the $k$-module
$$\Mag(V):= \bigoplus_{n\geqslant 1} k[Y_{n-1}]\otimes V^{\otimes n},$$
where $Y_{n-1}$ denote the set of planar binary trees with $n$
leaves. It is equipped with the following product
$$\big(t\otimes (x_1,\ldots ,\, x_m) \big)\diamond \big(s\otimes (x'_1,\ldots ,\, x'_n) \big):=
(t\vee s) \otimes (x_1,\ldots ,\, x_m,\,  x'_1,\ldots ,\, x'_n),$$
where $t$ belongs to $Y_{m-1}$ and $s$ belongs to $Y_{n-1}$. The
notation $t\vee s$ represents to the grafting of trees (\emph{cf.}
J.-L Loday \cite{Loday} Appendix A)
$$t\vee s = \vcenter{\xymatrix@R=8pt@C=8pt{t\ar@{-}[dr] &  & \ar@{-}[dl] s \\
 &*{}\ar@{-}[d] & \\
 & & }} $$

\begin{thm}[Free PostLie algebra]
\label{FreePostLie} The free PostLie algebra on the $k$-module $V$
is given by the module $\Li(\Mag(V))$, where the bracket $[,\, ]$
is the bracket coming from the free Lie algebra on $\Mag(V)$.
\end{thm}

\begin{deo}
The $\Postlie$ operad is a quadratic operad generated by two
operations $\circ$ and $[,\,]$. We have
$\Postlie=\F\big([,\,].sgn_{\Sy_2} \oplus \circ.
k[\Sy_2]\big)/(R)$, where the module of relations $R$ is the
direct sum $R_{[,\, ]}\oplus R_r \oplus R_l$. The module $R_{[,\,
]}$ corresponds to the Jacobi relation. The module $R_r$
corresponds to the first compatibility relation between $\circ$
and $[,\,]$, where the bracket is ``on the right'' and the module
$R_l$ corresponds to the second compatibility relation, where the
bracket is "on the left". The free operad
$\F\big([,\,].sgn_{\Sy_2} \oplus \circ. k[\Sy_2]\big)/(R)$ is
isomorphic the $\Sy$-module generated by binary trees where the
vertices are indexed by $\circ$ or $[,\,]$. The relation $R_r$
allows us to replace the following pattern
$$ 
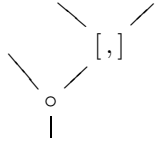
\xymatrix@R=8pt@C=8pt{ &\ar@{-}[rd] & &\ar@{-}[ld]  \\
{\  } \ar@{-}[rd]& & {[\, ,]}
\ar@{-}[dl]& \\
& {\circ} \ar@{-}[d]& & \\
& & & }$$ by a sum of trees with  vertices only indexed by
$\circ$. And the relation $R_l$ allows to replace the pattern
$$ 
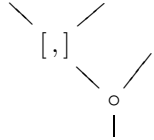
\xymatrix@R=8pt@C=8pt{ \ar@{-}[rd] & &\ar@{-}[ld] { \ }&  \\
& {[\, ,]} \ar@{-}[rd]& & \ar@{-}[dl] \\
& &{\circ} \ar@{-}[d]&  \\
& & & }$$ by a sum of trees where the vertex indexed by $\circ$ is
above the vertex indexed by $[\, ,]$. Therefore, the operad
$\Postlie$ is a quotient of the free module generated by binary
trees where the vertices are indexed by $\circ$ or $[\,,]$, such
that the vertices indexed by $\circ$ are above the vertices
indexed by $[\,,]$. With this presentation, the only remaining
relation is the Jacobi relation for the bracket $[\, ,]$. The
operad $\Postlie$ is isomorphic to the product $\Li \circ \Mag$,
which concludes the proof.$\cqfd$
\end{deo}

Once again, the $\Sy_n$-modules $\Postlie(n)$ are free
$k$-modules.

\begin{thm}
\label{ComtriKoszulPostlie}
The operad $\Postlie$ is the Koszul dual operad of the operad
$\Comtri$
\end{thm}

\begin{deo}
We have $\Comtri=\F(V)/(R)$ with $V=*.k[\Sy_2]\oplus \bullet.k$.
The module of relation is the direct sum $R=R_{\Pe}\oplus R_{\C}
\oplus R_{mix}$, where $R_{\Pe}$ denotes the Perm-relation of $*$,
$R_{\C}$ denotes the associativity of $\bullet$ and $R_{mix}$
represents the compatibility relations between $*$ and $\bullet$.
The operad $\Postlie$ is equal to $\F(W)/(S)$ with
$W=\circ.k[\Sy_2]\oplus [\,,].sgn_{\Sy_2}$. The module of
relations $S$ is the direct sum $S_{[\, ,]}\oplus S_{mix}$ where
$S_{[\, ,]}$ is the Jacobi relation and $S_{mix}$ represents the
compatibility relations.

The rank of the $k$-module $\F_{(2)}(V)=\F(V)(3)$ is $27$ and the
rank of $R$ is $9+2+9=20$ (for $R_{\Pe}$, $R_{\C}$ and $R_{mix}$).
We identify the linear dual of $*$ with $\circ$ and the linear
dual of $\bullet$ with $[\, ,]$. Therefore, we get $V^\vee\cong
W$. With this identification, one can easily see that $S\subset
R^\perp$. Since the rank of $S$ is $1+6=7=27-20$, the $k$-module
$S$ is equal to $R^\perp$ and $\Comtri^!=\F(V^\vee)/(R^\perp)\cong
\F(W)/(S)=\Postlie$. $\cqfd$
\end{deo}

$ $
\begin{flushleft}
3. $ $ {\bf  Koszul duality and homology of PostLie algebras}
\end{flushleft}
$ $
\begin{flushleft}
3.1. $ $ {\bf The Lie algebra associated to a PostLie algebra }
\end{flushleft}
$ $

\begin{pro}
Let $(L,\, \circ,\, [\, ,\, ])$ be a PostLie algebra. The bracket
defined by the formula
$$\{x,\, y\}:=x\circ y -y\circ x +[x,\,y] $$
is a Lie bracket. We denote by $L_{\{ , \}}$ the Lie algebra
$(L,\, \{ \, , \})$.
\end{pro}

\begin{deo}
The proof is given by a direct calculation. $\cqfd$
\end{deo}

\begin{pro}
A PostLie algebra $(L,\, \circ,\, [\, ,\, ])$ is a right module
over the Lie algebra $L_{\{ ,\,\}}$ for the following action
\begin{eqnarray*}
L\times L_{\{ ,\,\}} &\to& L \\
(x,\, y) &\mapsto& x\circ y.
\end{eqnarray*}
\end{pro}

\begin{deo}
Once again, the proof is given by a direct calculation. $\cqfd$
\end{deo}

One can extend the action of $L_{\{ ,\,\}}$ on $L$ to a right
action of the enveloping Lie algebra $\U(L_{\{ ,\,\}})$ by the
formula
\begin{eqnarray*}
L \times \U(L_{\{ ,\,\}})  &\to& L  \\
(g,\, v_1\otimes \cdots \otimes v_n) &\mapsto& ((g\circ v_1)
\cdots ) \circ v_n.
\end{eqnarray*}
We denote the action of $\U(L_{\{ ,\,\}})$ on $L$ by $\star$.

\begin{pro}
Let $L=(\Li(\Mag(V)),\, \circ,\, [,])$ the free PostLie algebra on
$V$. The module $\Mag(V) \subset L$ is a stable under the action
 of the Lie algebra $L_{\{ ,\,\}}$.
\end{pro}

\begin{deo}
It is an easy consequence of the avoiding pattern chosen in the
proof of Theorem \ref{FreePostLie}. $\cqfd$
\end{deo}

Therefore, $\Mag(V)$ is a right module over the enveloping algebra
$\U(L_{\{ ,\,\}})$.

\begin{pro}
\label{freeUmod} The right $\U(L_{\{ ,\,\}})$-module $\Mag(V)$ is
isomorphic to the free right $\U(L_{\{ ,\,\}})$-module generated
by $V$
$$\Mag(V) \cong V \otimes \U(L_{\{ ,\,\}}).$$
\end{pro}

\begin{deo}
We consider the following surjective morphism $\Phi$ of $\U(L_{\{
,\,\}})$-modules
\begin{eqnarray*}
V \otimes \U(L_{\{ ,\,\}}) &\to& \Mag(V) \\
v\otimes u &\mapsto& v\star u.
\end{eqnarray*}
We  give an inverse of $\Phi$. We define a magmatic structure on
$V\otimes\U(L_{\{ ,\,\}})$ by the product
$$(v\otimes u) \diamond (v'\otimes u'):= v \otimes (u\otimes i(v')\star u'), $$
where $i$ denotes the inclusion $V\hookrightarrow \Mag(V)$. Let
$j$ be the inclusion
\begin{eqnarray*}
V &\to& V \otimes \U(L_{\{ ,\,\}})\\
v &\mapsto& v\otimes 1.
\end{eqnarray*}
The morphism $\Phi$ preserves the binary products $\diamond$ and
$\circ$. We have
\begin{eqnarray*}
\Phi\big( (v \otimes u) \diamond (v'\otimes u')  \big) &=& v\star
(u\otimes (i(v')\star u')) = (v\star u) \star  (v'\star u')  \\
&=&  (v\star u) \circ  (v'\star u') = \Phi(v \otimes u) \circ
\Phi(v \otimes u').
\end{eqnarray*}
By definition of the free magmatic algebra, there exists a unique
morphism of magmatic algebras $\Psi$ such that the following
diagram commutes
$$\xymatrix{V \ar[r]^(0.4){i} \ar[rd]_{j}  & \Mag(V) \ar[d]^{\Psi}   \\
&  V \otimes \U(L_{\{ ,\,\}}). }$$

Since $\Phi \circ \Psi \circ i = i$, we have by definition of
$\Mag(V)$ that $\Phi \circ \Psi=id$.\\

It remains to prove that $\Psi$ is a morphism of right $\U(L_{\{
,\,\}})$-modules. It is enough to prove that $\Psi(t\star
u)=\Psi(t)\star u$ for $u\in L$. We denote by $\Li_{(k)}
(\Mag(V))$ the elements of the free PostLie algebra $L$ coming
from trees with $k$ vertices indexed by $[\, , ]$. We show the
previous statement by induction on $k$.

\begin{itemize}

\item[$\triangleright$] If $u$ belongs to
$\Li_{(0)}(\Mag(V))=\Mag(V)$, we have $\Psi(t\star u)=\Psi(t\circ
u)=\Psi(t)\diamond \Psi(u)$. We denote $\Psi(t)$ by $\sum_i
v_i\otimes u_i$ and $\Psi(u)$ by $\sum_j v'_i\otimes u'_j$.
Therefore, we have
\begin{eqnarray*}
\Psi(t)\diamond \Psi(u)&=& \sum_i v_i \otimes \big(u_i \otimes
\sum_j v'_j\star u'_j \big)=\sum_i v_i \otimes \big(u_i \otimes
\Phi \circ \Psi(u) \big) \\ &=& \sum_i v_i \otimes (u_i \otimes u
)=\Psi(t)\star u.
\end{eqnarray*}

\item[$\triangleright$] Suppose that the statement is true for
$k\leqslant n$. Let $u$ be an element of $\Li_{(n+1)}(\Mag(V))$.
In $\Li(\Mag(V))$, this element is a sum of elements of the form
$[w_1,\, w_2]$, where $w_1 \in \Li_{(k_1)}(\Mag(V))$ and $w_2 \in
\Li_{(k_2)}(\Mag(V))$ with $k_1, k_2 \leqslant n$. We have
\begin{eqnarray*}
\Psi(t\star u)&=& \Psi(t\star \sum [w_1,\, w_2])\\
&=& \sum \Psi(t\star (w_1\otimes w_2 -w_2\otimes w_1-w_1\circ
w_2+w_2\circ w_1)).
\end{eqnarray*}
Since $w_1\circ w_2$ and $w_2\circ w_1$ belong to
$\Li_{(m)}(\Mag(V))$ for $m\leqslant n$, we have for induction
hypothesis that
\begin{eqnarray*}
\Psi(t\star u)&=&  \sum \Psi(t)\star (w_1\otimes w_2 -w_2\otimes
w_1-w_1\circ w_2+w_2\circ w_1))\\ &=& \Psi(t)\star u.
\end{eqnarray*}
\end{itemize}

Since $\Psi\circ \Phi\circ j =j$, we have $\Psi\circ \Phi=id$ by
definition of the free $\U(L_{\{ ,\,\}})$-module on $V$, which
concludes the proof. $\cqfd$
\end{deo}

$ $
\begin{flushleft}
3.2. $ $ {\bf  Homology of a PostLie algebra}
\end{flushleft}
$ $

\begin{thm}[Homology of a PostLie algebra]
Let $(L,\, \circ,\, [\, ,\, ])$ be a PostLie algebra. Its operadic
homology theory is defined on the module
$\Comtri^\vee(L)=\bar{\Lambda}(V) \otimes \Lambda(V)$ by the
following boundary map :
\begin{eqnarray*}
&&d(x_1\wedge \ldots \wedge x_k\otimes y_1\wedge \ldots \wedge
y_l) = \\
&& \sum_{1\le i<j\le k} \pm \, x_1\wedge \ldots \wedge \lbrack
x_i,\, x_j \rbrack \wedge \ldots \wedge \widehat{x_j} \wedge
\ldots \wedge
x_k \otimes y_1\wedge \ldots \wedge y_l +\\
&& \sum_{1\le i<j\le k} \pm \, x_1\wedge \ldots \wedge
 x_k \otimes y_1\wedge \ldots \wedge \{ y_i,\, y_j \} \wedge \ldots \wedge \widehat{y_j}
\wedge \ldots\wedge y_l + \\
&& \sum_{1\le i \le k \atop 1\le j\le l} \pm \, x_1\wedge \ldots
 \wedge
\widehat{x_i} \wedge \ldots \wedge x_k \wedge (x_i\circ y_j)
\otimes y_1\wedge \ldots \wedge \widehat{y_j} \wedge \ldots \wedge
y_l,
\end{eqnarray*}
where the signs are given by the Koszul-Quillen rule.
\end{thm}

\begin{deo}
Since the Koszul dual of $\Postlie$ is the operad $\Comtri$, the
operadic homology for a PostLie algebra $L$ is given by the module
$\Postlie^{\ac}(L)=\Comtri^\vee(L)$ (\emph{cf.} section $2$). The
boundary maps are induced by the partial coproduct of the cooperad
$\Comtri^\vee$. The result is obtained by identifying the linear
bidual of $\circ$ to $*$ and the bidual of $\bullet$ to $[\, ,]$.
$\cqfd$
\end{deo}

We denote by $d_{[ ,]}(x_1\wedge \ldots \wedge x_k\otimes
y_1\wedge \ldots \wedge y_l)$ the first part $\sum_{1\le i<j\le k}
\pm \, x_1\wedge \ldots \wedge \lbrack x_i,\, x_j \rbrack \wedge
\ldots \wedge \widehat{x_j} \wedge \ldots \wedge x_k \otimes
y_1\wedge \ldots \wedge y_l$ of the previous boundary map and by
$d_{\{,\}}(x_1\wedge \ldots \wedge x_k\otimes y_1\wedge \ldots
\wedge y_l)$ the second part. Therefore, the map $d$ is equal to
$d_{[,]} +d_{\{,\}}$.

\begin{thm}[Homology of the free PostLie algebra]
Let $L=(\Li(\Mag(V)),\,  \circ, \, [,])$ the free Post-Lie algebra
on $V$. Its homology is equal to
$$H_l^{\Postlie}(L)=\left\{
\begin{array}{l}
V \quad \textrm{if} \quad l=1, \\
0 \quad \textrm{otherwise}.
\end{array}
\right. $$
\end{thm}

\begin{deo}
We consider the following filtration on the operadic complex
$(\bar{\Lambda}(V)\otimes \Lambda(V), d)$ of the free PostLie
algebra $L$
$$F_p:= \bar{\Lambda}(V)\otimes \Lambda_{\leqslant p}(V).$$
By the definition of the boundary map $d$, we get that this
filtration is preserved by $d$. We have $E^0_{p\,
q}=\Lambda_q(L)\otimes \Lambda_p(L)$ and $d_0=d_{[,]}$. Therefore,
the modules $E^1_{p,\, *}$ are equal to $H^{\Li}_*(\Li(\Mag(V)),\,
k)\otimes \Lambda_p(L)$, where $H^{\Li}_*(\Li(\Mag(V)),\, k)$ is
the Chevalley-Eilenberg homology of the free Lie algebra on
$\Mag(V)$ with coefficient in $k$. This homology is equal to
$$
E^1_{p,\, *}= \left\{
\begin{array}{ll}
\Mag(V)\otimes \Lambda_p(L) & \textrm{if} \quad *=1 \\
0 & \textrm{otherwise}.
\end{array} \right.
 $$
Since for $q=1$, we have $E^1_{p,\, q}=\Mag(V)\otimes
\Lambda_p(L)$ and $d_1=d{\{,\}}$, the second term of the spectral
sequence is isomorphic to
$$E^2_{p,\, 1}=H^{\Li}_p(L_{\{,\}}, \,
\Mag(V))=Tor_p^{\U(L_{\{ ,\,\}})}(\Mag(V),\, k).$$ We have seen in
the proposition~\ref{freeUmod} that $\Mag(V)$ is the free
$\U(L_{\{ ,\,\}})$-module on $V$. We finally get
$$E^2_{p,\, 1}=Tor_p^{\U(L_{\{ ,\,\}})}(\Mag(V),\, k)= \left\{
\begin{array}{ll}
V & \textrm{if} \quad p=0 \\
0 & \textrm{otherwise}.
\end{array}\right.$$
This spectral sequence is bounded. By the classical convergence
theorem of bounded spectral sequences, we have that $E^*$
converges to the homology of the free PostLie algebra on $V$ which
concludes the proof . $\cqfd$
\end{deo}

\begin{cor}
\label{ComtriKoszul}
 The operads $\Comtri$ and $\Postlie$ are
Koszul operads over any field $k$ of characteristic $0$.
\end{cor}

\begin{rem}
In \cite{CV}, we prove that the intervals of maximal length in the
posets of multi-pointed partitions are totally semi-modular. It
shows that they are Cohen-Macaulay over any field $k$ and over the
ring of integers $\mathbb{Z}$. As a consequence of
Theorem~\ref{KoszulCohen}, we get that the operads $\Comtri$ and
$\Postlie$ are Koszul over any field $k$ and over the ring of
integers $\mathbb{Z}$.
\end{rem}

\begin{center}
\textsc{Acknowledgements}
\end{center}
$ $

I wish to express my gratitude to Jean-Louis Loday for many useful
discussions on these subjects. I also would like to thank B.
Fresse for his help in my works and Roland Berger for pointing me
out the article of C. Cibils, P. Polo and W. Woodcock.\\

This paper was finished at the Mittag-Leffler Institut in
Stockholm, whose hospitality and financial support was greatly
appreciated.



{\small \textsc{Laboratoire J.-A. Dieudonn\'e, Universit\'e de
Nice Sophia-Antipolis, Parc Valrose, 06108 Nice Cedex, France}\\
E-mail address: \texttt{brunov@math.unice.fr}\\
URL: \texttt{http://math.unice.fr/$\sim$brunov}}
\end{document}